\documentclass[12pt,a4paper]{article}

\title{The Injective Spectrum of a Right Noetherian Ring I:\\Injective Spectra and Krull Dimension}
\author{Harry Gulliver}
\date{}

\usepackage[
top    = 2.5cm,
bottom = 2.5cm,
left   = 2.3cm,
right  = 2.3cm]{geometry}

\usepackage{upgreek}
\let\tau\uptau

\usepackage{amssymb}
\usepackage{amsmath}
\usepackage{amsthm}

\usepackage{tikz}
\usetikzlibrary{arrows}

\usepackage{fancyhdr}
\usepackage{enumerate}

\renewcommand{\qed}{\hfill$\blacksquare$\gap}
\newcommand{\gap}{\par \vspace{5mm}}

\let\take\smallsetminus

\newcommand{\im}{\mathrm{Im}}
\newcommand{\Mod}{\mathrm{Mod}}
\mathchardef\mhyphen="2D

\newcommand{\op}{\mathrm{op}}

\newcommand{\inj}{\mathrm{inj}}
\newcommand{\ann}{\mathrm{ann}}
\newcommand{\zg}{\mathrm{Zg}}
\newcommand{\spec}{\mathrm{Spec}}
\renewcommand{\mod}{\mathrm{mod}}
\newcommand{\ab}{\mathbf{Ab}}

\let\take\smallsetminus

\newcommand{\injspec}{\mathrm{InjSpec}}
\newcommand{\cl}{\mathrm{cl}}
\newcommand{\cd}{\mathrm{cd}}

\newcommand{\fp}{\mathrm{fp}}

\newtheorem{thm}{Theorem}[section]
\renewcommand{\proof}[1][]{\noindent\textsc{Proof:} \textit{#1}\par}

\newtheorem{lemma}[thm]{Lemma}
\newtheorem{cor}[thm]{Corollary}
\newtheorem{ex}[thm]{Example}
\newtheorem{prop}[thm]{Proposition}

\begin{document}

\maketitle
\pagenumbering{arabic}

\begin{abstract}
The injective spectrum is a topological space associated to a ring $R$, which agrees with the Zariski spectrum when $R$ is commutative noetherian. We consider injective spectra of right noetherian rings (and locally noetherian Grothendieck categories) and establish some basic topological results and a functoriality result, as well as links between the topology and the Krull dimension of the ring (in the sense of Gabriel and Rentschler). Finally, we use these results to compute a number of examples.
\end{abstract}

\tableofcontents

\section{Introduction and Background}

\subsection{Conventions}

Throughout, all rings will be associative and unital, but not necessarily commutative, and all modules will be unital right modules, unless otherwise specified. If $R$ is a ring, we denote by $\Mod\mhyphen R$ the category of all right $R$-modules, and by $\mod\mhyphen R$ the full subcategory of finitely presented modules. If $M$ is a module, or more generally an object in a Grothendieck category, we denote by $E(M)$ an injective hull of $M$. For modules (objects of a Grothendieck category) $L$ and $M$, we denote by $(L,M)$ the group of maps $L\to M$.

By ``functor'' we always mean ``additive, covariant functor.'' For a Grothendieck category $\mathcal{A}$, we denote by $\mathcal{A}^\fp$ the full subcategory of finitely presented objects; so $\mod\mhyphen R = (\Mod\mhyphen R)^\fp$.

\subsection{The Injective Spectrum}

Our starting point is the following classical result of Matlis:

\begin{thm}[\cite{matlis}, 3.1]\label{matlis bijection}
Let $R$ be a commutative noetherian ring. Then there is a bijection between prime ideals of $R$ and (isoclasses of) indecomposable injective modules, given by taking a prime $p$ to $E(R/p)$ in one direction and taking $E$ to the (unique) ideal which is maximal among annihilators of non-zero submodules of $E$, in the other.
\end{thm}

Building on this, Gabriel \cite[\S VI.3]{gabriel} defined a topology on the set of (isoclasses of) indecomposable injective modules over a ring $R$ which makes this bijection into a homeomorphism when $R$ is commutative noetherian and $\spec(R)$ is endowed with its Zariski topology. Gabriel's topology is defined as follows; given a module $M$, let $[M]$ denote the set of indecomposable injective modules $E$ such that $(M,E)=0$. Then a basis of open sets for the topology is given by the collection of all $[M]$ where $M\in\mod\mhyphen R$. We denote by $(M)$ the complement of $[M]$; \textit{i.e.}, the set of indecomposable injectives $E$ such that $(M,E)\neq 0$. See also \cite[\S 14.1.1]{prest} for a discussion of all this in more detail than in \cite{gabriel}.

We call the set of (isoclasses of) indecomposable injective $R$-modules together with this topology the (right) {\bf injective spectrum} of $R$ and denote it $\injspec(R)$. We refer to the topology as the {\bf Zariski topology}.

The Zariski spectrum of a commutative ring carries a sheaf of rings; we explore sheaves of rings on the injective spectrum in the sequel paper \cite{gulliver2}, but for now consider the Zariski and injective spectra purely as topological spaces.

Of course, the description of the injective spectrum does not depend in any way on the ring itself, only on $\Mod\mhyphen R$ (and so, in particular, Morita equivalent rings have homeomorphic injective spectra); as such, given any Grothendieck category $\mathcal{A}$, we can define the injective spectrum $\injspec(\mathcal{A})$ to be the set of (isoclasses of) indecomposable injective objects of $A$, with the topology having for basis of open sets the collection of all $[A]=\{E\in\injspec(\mathcal{A})\mid (A,E)=0\}$ as $A$ ranges over $\mathcal{A}^\fp$, the category of finitely presented objects of $\mathcal{A}$. Note that there is an abuse of notation here: technically, $\injspec(R)$ is actually $\injspec(\Mod\mhyphen R)$.

 Recall that a Grothendieck category is {\bf locally noetherian} if there is a generating set of noetherian objects; for a ring $R$, $\Mod\mhyphen R$ is locally noetherian if and only if $R$ is right noetherian. Although our principal interest here is in injective spectra of right noetherian rings, many of our results hold for general locally noetherian Grothendieck categories, so we prove them in this generality.

There is an alternative topology on the injective spectrum. The full story is that there is a topological space associated to a ring (or suitable Grothendieck category) called the Ziegler spectrum, introduced in \cite{ziegler}; the injective spectrum is a subset of the Ziegler spectrum, and when the Grothendieck category is locally noetherian, there is a useful relationship between the Zariski topology on the injective spectrum and the restriction of the Ziegler topology. We do not require the Ziegler spectrum in full, so simply present some results in the generality relevant to us.

Let $\mathcal{A}$ be a locally noetherian Grothendieck category. The {\bf Ziegler topology} on the injective spectrum $\injspec(\mathcal{A})$ is the topology having $\{(A)\mid A\in\mathcal{A}^\fp\}$ as a basis of open sets (recall that $(A)=\{E\in\injspec(\mathcal{A})\mid (A,E)\neq 0\}$). So the basic open sets of the Ziegler topology are the complements of the basic open sets for the Zariski topology; for this reason, the Zariski topology on the injective spectrum is sometimes called the {\bf dual-Ziegler topology}. The definition of the Ziegler topology given here is not the standard one, but is equivalent by \cite[5.1]{prest_rajani}.

A key fact we shall need is the following:

\begin{lemma}\label{ziegler compact opens}
Let $\mathcal{A}$ be a locally noetherian Grothendieck category. Then the basic open sets $(A)$ for $A\in\mathcal{A}^\fp$ on $\injspec(\mathcal{A})$ are compact open in the Ziegler topology.
\end{lemma}

\proof
We give a sketch proof only, which requires knowledge of the ZIegler spectrum; references are included to assemble a complete proof.

In the full Ziegler spectrum, the basic open sets given by pp-pairs are precisely the compact opens, by \cite[Thm. 4.9]{ziegler}. By \cite[5.1]{prest_rajani}, these restrict to the given basic open sets on the injective spectrum, for $\mathcal{A}$ locally noetherian. Moreover, in this case the injective spectrum is a closed subset of the full Ziegler spectrum (see \textit{e.g.}, \cite[\S 5.1.1]{prest}). The result follows.\qed

There is a further model-theoretic result we require. Given a ring $R$ and a right $R$-module $M$, say that the system
\[\left(\sum_i x_ir_{ij}=m_j\right)_j\]
 of $R$-linear equations in $M$ is {\bf consistent} if whenever a finite linear combination of the left-hand sides is 0, the corresponding linear combination of the right-hand sides is also 0. That is, for any collection of $s_j\in R$ almost all zero, we have
\[\forall i \left( \sum_j r_{ij}s_j = 0 \right) \quad\Rightarrow\quad \sum_j m_js_j=0.\]

We have the following result by Eklof and Sabbagh.

\begin{thm}[\cite{eklof}, \S3]\label{equations in injectives}
An $R$-module $M$ is injective if and only if every consistent system of equations in $M$ has a solution in $M$. By Baer's Criterion, it is sufficient to consider systems of equations in just one variable.
\end{thm}

\subsection{Pappacena's ``Weak Zariski'' Topology}\label{weak zariski}

In \cite{pappacena}, Pappacena studies a topological space he calls the injective spectrum of a noncommutative space. This is closely related to the notion of injective spectrum we consider here, both in provenance and details, but is not the same. For the avoidance of confusion, we clarify here the distinction.

Pappacena defines a {\bf noncommutative space} to be a Grothendieck abelian category, and focuses on those which are locally noetherian. A {\bf weakly closed subspace} is then defined to be a full subcategory which is itself Grothendieck, and is closed under subquotients, coproducts, and isomorphisms, and is such that the inclusion functor admits a right adjoint.

Given a noncommutative space $X$, Pappacena then defines the injective spectrum to be the set of isoclasses of indecomposable injective objects in $X$, with the {\bf weak Zariski topology} defined as follows. For each weakly closed subspace $Z$, we define $V(Z)$ to be the set of indecomposable injectives containing a subobject from $Z$. Then the collection of all sets $V(Z)$ for $Z$ weakly closed is a basis of closed sets for the weak Zariski topology. The space with this topology is what Pappacena refers to as the injective spectrum.

Clearly the points of Pappacena's injective spectrum coincide with ours. Moreover, for any module $M$, there is a weakly closed subspace $\sigma[M]$ (see \cite[\S15]{wisbauer}), which consists of all subquotients of coproducts of copies of $M$. Then the weak Zariski basic closed set $V(\sigma[M])$ is precisely the set $(M)=\{E\mid (M,E)\neq 0\}$. In particular, taking $M$ to be finitely presented, every basic closed set in our topology is still basic closed for Pappacena; so, despite the name, the weak Zariski topology refines the Zariski topology.

This refinement is however strict, in general, because Pappacena does not require any finite presentation condition on his basic open sets. To see that the refinement is strict, work in the category of abelian groups, since $\injspec(\mathbb{Z})\cong\spec(\mathbb{Z})$, and take $M$ to be the direct sum of each simple group $\mathbb{Z}/(p)$ for $p$ an odd prime. Then $V(\sigma[M])$ consists of all indecomposable injective abelian groups except $\mathbb{Q}$ and $\mathbb{Z}_{2^\infty}$. To prove this, show that $\mathbb{Q}$ and $\mathbb{Z}/(2)$ are torsionfree for the hereditary torsion class generated by $M$, so $(M,\mathbb{Q})=0=(M,\mathbb{Z}_{2^\infty})$. But this set is not Zariski-closed, for it contains the generic point $\mathbb{Q}$ but not the whole space.

Pappacena also mentions a second topology on the injective spectrum, which he calls the {\bf strong Zariski topology}. However, he points out that this topology is trivial for simple rings (remarks after Prop. 4.11 of \cite{pappacena}), in contrast with the Zariski topology studied here (see Theorem \ref{spectrum of a 1-critical ring} and following remarks), so this topology is also different to ours.

\subsection{Outline of Paper}

This paper is essentially one half of the author's PhD thesis, which was prepared under the supervision of Prof Mike Prest and submitted to the University of Manchester in June 2019. The other half forms a second paper, \cite{gulliver2}, covering sheaves of rings and modules over the injective spectrum, links with hereditary torsion theories, and further topological results.

In section \ref{prelim results}, we establish a criterion for specialisation within the injective spectrum and consider closed points. We also prove that certain ring maps $R\to S$ induce continuous maps $\injspec(S)\to \injspec(R)$, though we cannot hope to have such a result for all ring maps.

In section \ref{krull dim}, we consider Krull dimension (in the sense of Gabriel \& Rentschler) and its relation to the topological structure of the injective spectrum. We use this dimension to prove results regarding specialisation, points of maximal dimension, and irreducibility of the injective spectrum of a noetherian domain.

Finally, in section \ref{first examples}, we use the machinery built up to compute a number of examples. We determine in detail the injective spectra of all right artinian rings, all right noetherian domains of Krull dimension 1, and (the universal enveloping algebra of) the first Heisenberg algebra. We also exhibit more complicated behaviour in the injective spectrum of the quantum plane, showing that certain conjectures one might reasonably form are in fact false.

The sequel paper \cite{gulliver2} builds on this work, but in a different direction; it considers sheaves of rings and modules on the injective spectrum, and relates the injective spectrum to hereditary torsion theories, exploiting this link to prove further results about the topology on the spectrum.

\subsection{Acknowledgements}

I owe an enormous debt of gratitude to Prof Mike Prest, who supervised my PhD, in which I completed the work which forms this paper; his guidance and patience have been exemplary. Thanks are also due to Tommy Kucera, Ryo Kanda, Lorna Gregory, and Marcus Tressl, as well as my examiners, Omar Le{\'o}n S{\'a}nchez and Gwyn Bellamy, for stimulating discussions and helpful feedback. Finally, I am grateful to EPSRC for providing the funding for my PhD.

\section{Preliminary Results}\label{prelim results}

\subsection{Specialisation; Closed Points}

Throughout this section, we fix a locally noetherian Grothendieck category $\mathcal{A}$. The reader may take $\mathcal{A}$ to be the category of right modules over a right noetherian ring, if preferred.

First we consider the closure of a point in the injective spectrum. For $E,F\in\injspec(\mathcal{A})$, we write $E\leadsto F$ and say that $E$ {\bf specialises to} $F$ if $F\in \cl(E)$ - \textit{i.e.}, if every closed set containing $E$ also contains $F$.\gap

\begin{lemma}\label{specialisation tfae}
For $E,F\in \injspec(\mathcal{A})$, the following are equivalent:
\begin{enumerate}
\item $E\leadsto F$;
\item For every finitely presented object $A$, if $(A,E)\neq 0$, then $(A,F)\neq 0$;
\item For every object $A$, if $(A,E)\neq 0$, then $(A,F)\neq 0$;
\item There is a cardinal $\lambda$ such that $E$ embeds in $F^\lambda$.
\end{enumerate}
\end{lemma}

\proof
\noindent $(1\Leftrightarrow 2)$:

Since the sets $(A)$ for $A$ finitely presented form a basis of closed sets, $E\leadsto F$ if and only if for every $A$ finitely presented, $E\in (A)$ implies $F\in(A)$; but this is precisely statement (2).

\noindent $(2\Rightarrow 3)$:

Suppose $A$ is any object of $\mathcal{A}$ and $(A,E)\neq 0$. Take $f:A\to E$ a non-zero morphism; then there is a noetherian generating object $G$ and a map $g:G\to A$ such that $f\circ g$ is non-zero. Then the image $B$ of $g$ is a finitely presented object which maps to $E$, so by (2) we have $(B,F)\neq 0$. Since $F$ is injective, any non-zero morphism $B\to F$ extends to a morphism $A\to F$, establishing (3).

\noindent $(3\Rightarrow 4)$:

Take $\lambda$ to be the cardinality of $(E,F)$; then the product of all morphisms $E\to F$ gives a morphism $\phi:E\to F^\lambda$. We show that $\phi$ is a monomorphism; to do this, we fix a generator $G$ for $\mathcal{A}$ and show that for every non-zero morphism $\psi:G\to E$, we have $\phi\circ\psi\neq 0$; this suffices.

So suppose that $\psi$ is a non-zero morphism $G\to E$; then $(\psi(G),E)\neq 0$, so $(\psi(G),F)\neq 0$, by (3). Since $\psi(G)$ is a subobject of $E$ and $F$ is injective, we can extend any non-zero morphism $\rho:\psi(G)\to F$ to a morphism $\hat{\rho}:E\to F$ which does not vanish on $\psi(G)$. Therefore $\hat{\rho}\circ\psi\neq 0$; but $\phi$ is the product of all morphims $E\to F$, so $\hat{\rho}$ is a component of $\phi$ and so $\phi\circ\psi\neq 0$, as required.

\noindent $(4\Rightarrow 2)$:

Suppose $A$ is finitely presented and $(A,E)\neq 0$. Then we can take $\phi:A\to E$ non-zero and $\phi(A)$ embeds in $E$; hence $\phi(A)$ embeds in some power $F^\lambda$, by $4$. But then $(\phi(A),F)\neq 0$, so $(A,F)\neq 0$.\qed

\begin{cor}\label{specialisation lemma}
Suppose $A$ is a uniform object of $\mathcal{A}$ (so that $E(A)$ is indecomposable) and $E\in\injspec(\mathcal{A})$. If $A$ embeds in a product of copies of $E$, then $E(A)\leadsto E$.
\end{cor}

\proof
Take an embedding $\phi:A\hookrightarrow E^\lambda$, for some cardinal $\lambda$. Since $E$ is injective, this extends to a morphism $\hat{\phi}:E(A)\to E^\lambda$, which restricts to $\phi$ on $A$. Since $\ker(\phi)=0$, we have $\ker(\hat{\phi})\cap A=0$; but $A$ is essential, so $\hat{\phi}$ is mono. Then part (4) of the above Lemma \ref{specialisation tfae} gives the conclusion.\qed

This Corollary will be our principal tool for establishing specialisation in examples. Note that statement (4) of Lemma \ref{specialisation tfae} says that $E$ is in the hereditary torsionfree class cogenerated by $F$, and the Corollary says that it suffices to check that any non-zero submodule is torsionfree. We do not consider torsion theories or their relation to the injective spectrum in this paper, but they form a principal theme of the sequel paper \cite{gulliver2}.\gap

For $R$ commutative noetherian, the closed points of $\spec(R)$ are the maximal ideals, so the closed points of $\injspec(R)$ are the points of the form $E(R/m)$ for $m$ a maximal ideal. Of course, $m$ is maximal if and only if $R/m$ is a simple module, so the closed points are those which contain a simple submodule. In general, we have the following:

\begin{prop}\label{socle implies closed}
Let $E\in\injspec(\mathcal{A})$ and suppose that $E$ has a simple subobject $S$. Then $E$ is a closed point.\end{prop}

\proof
Suppose that $E\leadsto F$; we show that $E\cong F$. The embedding $S\to E$ gives a non-zero element of $(S,E)$, so $(S,F)$ must be non-zero, by part (3) of Lemma \ref{specialisation tfae}. Since $S$ is simple, this implies that $S$ embeds in $F$; but $F$ is indecomposable injective, so $F$ is an injective hull of $S$, and so too is $E$. So $E\cong F$.\qed

We can also obtain a partial converse to this.

\begin{prop}\label{basic closed implies socle}
Let $E$ be a closed point in $\injspec(\mathcal{A})$ such that $\{E\}=(A)$ is basic closed (for some finitely presented object $A$). Then $E$ has a simple subobject. Moreover, if $\injspec(\mathcal{A})$ is a noetherian space, then the additional hypothesis that $\{E\}$ be basic closed can be dropped.
\end{prop}

\proof
Since $A$ is finitely generated in a locally noetherian category, it has at least one simple quotient. For any simple quotient $S$ of $A$, we have $(A,E(S))\neq 0$, so $E(S)\in (A)=\{E\}$, so $E=E(S)$. Therefore $E$ has a simple subobject, as required.

Now suppose that $\injspec(\mathcal{A})$ is noetherian. Then, since $E$ is a closed point, $\{E\}$ is an intersection of basic closed sets, and this intersection can be taken to be finite, by noetherianity. Now, each basic Zariski-closed set $(A)$ is Ziegler-open, and so $\{E\}$ must be Ziegler-open. Hence $\{E\}$ can be written as a union of basic Ziegler-open sets $(A)$. But, since $\{E\}$ is a singleton, it must in fact be equal to a single basic Ziegler-open set, which is the same as a basic Zariski-closed set. So in this case the hypothesis that $\{E\}$ be basic open is automatic.\qed

We will see in section \ref{first examples} that the injective spectrum of a right noetherian domain of Krull dimension 0 or 1 is noetherian (see section \ref{krull dimension} for Krull dimension). However, subsection \ref{quantum plane} exhibits a right noetherian domain of Krull dimension 2, whose injective spectrum fails to be noetherian; moreover, there is a closed point which contains no simple submodule, in contrast with the above Proposition.\gap

\subsection{Functoriality}

Of course, the Zariski spectrum is not simply a topological space associated to a commutative ring. It is a contravariant functor from the category of commutative rings to the category of topological spaces.

Unfortunately, there are serious limits to what we can hope for when it comes to functoriality of $\injspec$. In particular, we cannot hope for a duality of categories, as the following example, shown to me by Ryo Kanda, illustrates.

\begin{ex}\label{kanda's matrix example}
Let $k$ be a field and $M_2(k)$ denote the ring of $2\times 2$ matrices with coefficients in $k$. Consider the ring map $f:k\times k\to M_2(k)$, embedding the leading diagonal. We show that $f$ cannot induce a natural map of topological spaces $\injspec(M_2(k))\to\injspec(k\times k)$.
\end{ex}
Since $M_2(k)$ is simple artinian, every module is injective. There is a unique indecomposable module, $k^2$, so $\injspec(M_2(k))=\{k^2\}$. On the other hand, indecomposable $k\times k$-modules are just indecomposable modules over each factor, and again all modules are injective, since the ring is semisimple artinian, so $\injspec(k\times k)$ is a 2-point, discrete space.

So any map induced by $f$ on injective spectra would need to take the single point of $\injspec(M_2(k))$ to one of the two points of $\injspec(k\times k)$. But there is nothing to distinguish these points; indeed, there is an automorphism of $k\times k$ swapping the two factors, so there is no natural way to pick one out to be the image of an induced map.

In particular, if there were a duality of categories between the category of rings and some subcategory of topological spaces, as in the commutative case, then the monomorphism $f$ would induce a topological epimorphism. But there is no epimorphism from a one point space to a discrete, two point space.\qed

Although we cannot have functoriality in general, we do get the following results. By a {\bf topological embedding}, we mean a continuous map which is a homeomorphism onto its image.

\begin{thm}\label{induced continuous map}
Let $\mathcal{A}$ and $\mathcal{B}$ be Grothendieck abelian categories and let $F:\mathcal{A}\to\mathcal{B}$ be a fully faithful functor with an exact left adjoint $L:\mathcal{B}\to \mathcal{A}$ which preserves finitely presented objects. Then $F$ induces by restriction and corestriction a continuous injection $\injspec(\mathcal{A})\to\injspec(\mathcal{B})$.

If, moreover, $L$ is such that for any finitely presented object $A$ of $\mathcal{A}$ there is $B\in\mathcal{B}$ finitely presented with $L(B)=A$, then the induced map is a topological embedding.

Note that if $\mathcal{B}$ is locally finitely presented, $L$ is the localisation at a hereditary torsion theory of finite type, and $F$ is the adjoint inclusion, then all of these assumptions are satisfied.\end{thm}

\proof
Since $F$ is fully faithful, it preserves indecomposables. For if $A\in\mathcal{A}$ is such that $F(A)$ is decomposable, then there is a non-trivial idempotent endomorphism $e$ of $F(A)$ (projection onto a summand). Then, since $F$ is full, $e=F(\epsilon)$ for some $\epsilon:A\to A$, with $F(\epsilon^2-\epsilon)=e^2-e=0$; but $F$ is faithful, so $\epsilon^2-\epsilon=0$, so $\epsilon$ is an idempotent, and necessarily non-trivial. So $A$ is decomposable.

Since $F$ has an exact left adjoint, it preserves injectives. For if $E\in\mathcal{A}$ is an injective, then $(-,F(E))=(L(-),E)$ is the composition of the exact functors $(-,E)$ and $L$, so is exact, implying that $F(E)$ is injective.

So $F$ gives a well-defined function $\injspec(\mathcal{A})\to\injspec(\mathcal{B})$, which is injective since $F$ is fully faithful, hence injective on isoclasses of objects.

To show that this map is continuous, we show that the preimage of any basic open set $[B]$, for $B\in\mathcal{B}^\mathrm{fp}$, is open. To see this, note that $F^{-1}[B]=\{E\in\injspec(\mathcal{A})\mid (B,F(E))=0\}$; but $(B,F(E))=(L(B),E)$, so $F^{-1}[B]=[L(B)]$. Since $L$ preserves finitely presented objects, this is open (even basic open).

Finally, note that for any $A\in\mathcal{A}^\mathrm{fp}$, $F[A]=\{F(E)\mid E\in\injspec(\mathcal{A}) \wedge (A,E)=0\}$. If $A=L(B)$, for some $B\in\mathcal{B}^\mathrm{fp}$, then $(A,E)=(L(B),E)=(B,F(E))$, so then $F[A]=\{F(E)\mid E\in \injspec(\mathcal{A})\wedge (B,F(E))=0\}=[B]\cap F(\injspec(\mathcal{A}))$. So if for any $A$ we can find such a $B$, then the image of any basic open set is relatively open within the image and so $F$ is an embedding of topological spaces.\qed

\begin{cor}\label{induced map}
Let $R$ and $S$ be rings and $f:R\to S$ a ring epimorphism such that $S$ is flat when viewed as a left (\textit{sic}) $R$-module. Then $f$ induces a continuous injection of topological spaces $f^\ast:\injspec(S)\to\injspec(R)$ by restriction of scalars. If, moreover, every finitely presented $S$-module can be obtained by extending scalars on a finitely presented $R$-module, then $f^\ast$ is a topological embedding.\end{cor}

\proof
The restriction of scalars functor $\Mod\mhyphen S\to \Mod\mhyphen R$ is faithful and, since $f$ is epi, is also full \cite[Corollary 1.3]{silver}. It has a left adjoint $-\otimes_R S$, which is exact, since ${}_RS$ is flat. Therefore the above Theorem applies.\qed

Note that any Ore localisation fulfills the hypotheses of this Corollary.

We can relax the requirement that ${}_R S$ be flat, but at the cost that we must strengthen $f$ from an epimorphism to a surjection with centrally-generated kernel. We show this below, but first require a Lemma.

\begin{lemma}\label{central elements act nilpotently}
Let $R$ be a right noetherian ring, $I$ a two-sided ideal of $R$ and $F$ an indecomposable injective $R/I$-module. Suppose that $z\in I\cap Z(R)$ is a central element of $I$. Then for any $x\in E(F_R)$, there is a natural number $n$ such that $xz^n=0$.
\end{lemma}

\proof
This proof is adapted from an argument of E. Noether, as presented in Theorem 3.78 of \cite{lam}.

If $x=0$, there's nothing to prove; so suppose $x\neq 0$ and let $Q:=\ann_R(x)$. Then $Q$ is a proper right ideal of $R$ and $xR\cong R/Q$. Since $F$ is uniform, so is $E(F_R)$, and therefore so is $xR$, so $Q$ is $\cap$-irreducible among right ideals of $R$.

Also by uniformity, $xR\cap F_R\neq 0$; so let $D$ be the right ideal of $R$ containing $Q$ such that $D/Q$ corresponds to $xR\cap F_R$ under the isomorphism $R/Q\cong xR$. Since $F_RI=0$, $(D/Q)I=0$, so $DI\subseteq Q$; in particular, $Dz\subseteq Q$. Moreover, $D$ strictly contains $Q$, since $xR\cap F_R\neq 0$.

Now consider the sets $(Q:z^i):=\{r\in R\mid rz^i\in Q\}$. These are right ideals, since $z$ is central, and $(Q:z^i)\subseteq (Q:z^{i+1})$ for all $i$. Since $R$ is right noetherian, this chain stabilises, and so there is some $n$ such that $(Q:z^n)=(Q:z^{n+1})$. We show that, for this particular $n$,
\[Q=(Q+z^nR)\cap D.\]

The left-to-right inclusion is clear. Let $y$ be an element of the right-hand side; so for some $q\in Q$, and $r\in R$, $y=q+rz^n$ and $y\in D$. Then $yz\in Q$, since $Dz\subseteq Q$. So $yz=qz+rz^{n+1}\in Q$, hence $rz^{n+1}=yz-qz\in Q$. So $r\in(Q:z^{n+1})=(Q:z^n)$, so $rz^n\in Q$. But then $y=q+rz^n\in Q$, proving the claim.

Now, since $D$ strictly contains $Q$ and $Q$ is $\cap$-irreducible, we must have $Q=Q+z^nR$, so $z^n\in Q=\ann_R(x)$, so $xz^n=0$.\qed

\begin{thm}\label{induced by surjection}
Let $R$ and $S$ be right noetherian rings and let $f:R\to S$ be a surjection, with $\ker(f)$ generated by central elements. Then $f$ induces a topological embedding $f^\ast:\injspec(S)\to \injspec(R)$ given by $F\mapsto E(F_R)$, where $E$ denotes the injective hull (taken, in this case, in $\Mod\mhyphen R$).
\end{thm}

\proof
Since $f$ is surjective, the submodule structure of any $S$-module is the same when scalars are restricted along $f$ to $R$. In particular, for any uniform $S$-module $U$, $U_R$ is uniform. Therefore, for $F\in\injspec(R)$, $E(F_R)$ is indecomposable. So $f^\ast$ is well-defined.

Moreover, if $F,G\in \injspec(S)$ and $E(F_R)=E(G_R)$, then $F_R$ and $G_R$ have a common submodule (up to isomorphism), by uniformity. So there are $M_R\leq F_R$, $N_R\leq G_R$ with $M_R\cong N_R$. But $F$ and $G$ have the same submodules over $R$ and over $S$, so $M_R$ and $N_R$ are indeed restrictions of some $S$-submodules $M_S$ and $N_S$, and the isomorphism between them is also an isomorphism of $S$-modules, since $f$ is onto. So $F_S=E(M_S)$ and $G_S=E(N_S)$ are isomorphic and so $f^\ast$ is injective.\gap

Next we show that $f^\ast$ maps closed sets to relatively closed sets. So let $M\in\mod\mhyphen S$, so that $(M)$ is a basic closed set in $\injspec(S)$. We show that $f^\ast(M)=(M_R)\cap \im(f^\ast)$. Since $f$ is a surjection, a finite generating set for $M$ as an $S$-module will still generate it over $R$, and $R$ is right noetherian, so $M_R$ is finitely presented and hence $(M_R)$ is (basic) closed.

On the one hand, if $E(F_R)\in f^\ast(M)$, then $(M,F)\neq 0$, so $(M_R,F_R)\neq 0$ (since restriction of scalars is faithful) and so $(M_R,E(F_R))\neq 0$ and so $E(F_R)\in (M_R)\cap\im(f^\ast)$.

On the other hand, if $E(F_R)\in (M_R)\cap \im(f^\ast)$, then $(M_R,E(F_R))\neq 0$, so $M_R$ has a submodule $N_R$ such that $(N_R,F_R)\neq 0$. Since $M$ has the same submodule structure over $R$ and over $S$, $M$ has an $S$-submodule $N$ restricting to $N_R$ and any non-zero map $N_R\to F_R$ is also $S$-linear, so $(N,F)\neq 0$. Then injectivity of $F$ allows us to extend any non-zero, $S$-linear map $N\to F$ to a map $M\to F$, showing that $F\in (M)$ and so $E(F_R)\in f^\ast(M)$.

Finally, since $f^\ast$ is injective, for any sets $A_i$, $f^\ast\left(\bigcap A_i\right)=\bigcap f^\ast (A_i)$, so this extends from basic closed sets to arbitrary closed sets.

Note that, up to this point, we have not used our assumption on the kernel.\gap

Finally, we assume that $\ker(f)$ is generated by central elements and prove that $f^\ast$ is continuous. It suffices to prove when $\ker(f)$ is generated by a single central element; for the quotient by any $n$ central elements factors through quotienting out a single central element at a time, $n$ times. So suppose $\ker(f)=zR$ for some $z\in Z(R)$.

It suffices to prove that for any finitely presented $R$-module $M$, the preimage of the basic closed set $(M)$ under $f^*$ is basic closed. The preimage is $(f^*)^{-1}((M))=\{F\in\injspec(S)\mid (M,E(F_R))\neq 0\}$. We show that $(f^*)^{-1}((M))=(M\otimes_R S)$, which is basic closed in $\injspec(S)$, since tensoring preserves finitely presented modules.

Certainly if $F\in (M\otimes_R S)$, then $(M,F_R)=(M\otimes_RS,F)\neq 0$, so after embedding in $E(F_R)$, we see that $(M,E(F_R))\neq 0$, so $F\in (f^*)^{-1}((M))$.

Conversely, if $(M,E(F_R))\neq 0$, take $\phi:M\to E(F_R)$ a non-zero map and let $m_1,\hdots,m_n$ be a generating set for $M$. By Lemma \ref{central elements act nilpotently}, for each $i$ there is some smallest integer $\nu_i\geq 0$ such that $\phi(m_i)z^{\nu_i}=0$. Let $\nu=\max_{i=1}^n\{\nu_i\}$. We cannot have all $\nu_i$ zero, since then $\phi$ would be zero, a contradiction; so $\nu\geq 1$.

So $\phi(M)z^\nu=0$ and $\phi(M)z^{\nu-1}\neq 0$; so $\phi(M)z^{\nu-1}$ is a non-zero submodule of $A:=\ann_{E(F_R)}(z)$. We show that $A= F_R$. Certainly $F_R\subseteq A$, since $F$ is a module over $S\cong R/zR$, so $F_Rz=0$. Moreover, $Az=0$, so $A$ naturally has an $S$-module structure. But $F_R\leq A\leq E(F_R)$, so $A$ is an essential extension of $F_R$; since $f:R\to S$ is surjective, the submodule structure of $A$ does not depend on whether we regard it as an $R$-module or $S$-module, so $A_S$ is an essential extension of $F$. But $F$ is injective, so admits no proper essential extension as an $S$-module; therefore $A=F$, as claimed.

So $\phi(M)z^{\nu-1}$ is a non-zero submodule of $F_R$. Let $\psi:E(F_R)\to E(F_R)$ be the ``multiplication by $z^{\nu-1}$'' map, which is $R$-linear, since $z$ is central. Then $\psi\circ\phi:M\to E(F_R)$ is non-zero and has image in $F_R$, so $(M,F_R)\neq 0$, and therefore $F\in (M\otimes_R S)$, completing the proof.\qed

\section{Krull Dimension}\label{krull dimension}

\subsection{Definitions of Krull Dimension and Critical Dimension; Basic Results}

There is a strong relationship between the Krull dimension of a commutative ring (and its quotients) and the topological dimension (in the sense of chains of irreducible closed sets) of its Zariski spectrum and subsets thereof. We therefore seek a notion of dimension for noncommutative rings which generalises this to the injective spectrum.

As in the commutative case, we shall work with Krull dimension; however, we use the ``noncommutative Krull dimension'' developed by Gabriel and Rentschler \cite{gabrentsch}, followed by Gordon and Robson \cite{gr}. We recall the definitions and basic results.

We first define the {\bf deviation} of a poset (from being artinian), by a transfinite induction. Let $P$ be a poset.

For $a,b\in P$, let $[a,b]$ denote the subposet of $P$ defined by the formula $a\leq x\leq b$. If $a_0>a_1>\hdots$ is a chain in $P$ of type $\omega^\op$, we call the subposets $[a_{i+1},a_i]$ for $i<\omega$ the {\bf factors} of the chain.

If $P$ is trivial (has no comparable elements), we say $P$ has deviation $-1$ (some authors prefer $-\infty$). If $P$ is non-trivial and has the descending chain condition, then we say it has deviation 0; otherwise it has deviation at least 1. Having defined what it means to have deviation at least $\alpha$, we say that $P$ has deviation at least $\alpha+1$ if there is an infinite descending chain in $P$ of type $\omega^\op$ all of whose factors have deviation at least $\alpha$. If $\lambda$ is a limit ordinal and $P$ has deviation at least $\alpha$ for all $\alpha<\lambda$, we say that $P$ has deviation at least $\lambda$. Finally, $P$ has deviation $\alpha$ if it has deviation at least $\alpha$ and does not have deviation at least $\alpha+1$.

Note that some posets do not have well-defined deviation. For instance, suppose $P=(\mathbb{Q}\cap [0,1],<)$ has deviation $\alpha$, where here $[0,1]$ denotes the closed unit interval in $\mathbb{R}$. Then take any chain $a_0>a_1>\hdots$; each factor $[a_{i+1},a_i]$ is isomorphic to $P$, so has deviation $\alpha$; but then $P$ has deviation at least $\alpha+1$, a contradiction. Indeed, a poset fails to have deviation if and only if it contains a dense linear order \cite[Proposition 6.1.12]{mcr}.

For $A$ an object in a Grothendieck category, we define the {\bf Krull dimension} $K(A)$ to be the deviation of the poset of subobjects of $A$. So $K(A)$ measures how far $A$ is from being artinian. We define the (right) Krull dimension of a ring to be its dimension as a (right) module over itself: $K(R):=K(R_R)$.

For sufficiently ``large'' objects, whose poset of subobjects contains a dense linear order, the Krull dimension need not be defined \cite[Lemma 6.2.6]{mcr}. In the case where $K(A)$ fails to exist, we write $K(A)=\infty$. For the purposes of inequalities, we regard $\infty$ as being strictly greater than any ordinal.

Let $\alpha$ be an ordinal. We say that a non-zero object $A$ is $\pmb{\alpha}${\bf -critical} if $K(A)=\alpha$ and for any non-zero subobject $B$ of $A$, $K(A/B)<\alpha$. We say that $A$ is {\bf critical} if there exists some $\alpha$ such that $A$ is $\alpha$-critical. Note that a $0$-critical object is precisely a simple object.

Some basic facts about Krull dimension are summarised below.

\begin{prop}[\cite{mcr}, \S\S 6.1-6.3]\label{krull dim}
Let $\mathcal{A}$ be a Grothendieck category with generator $G$, and $A$ an object of $\mathcal{A}$.
\begin{enumerate}
\item If $A$ is noetherian, then $K(A)$ exists.
\item If $B$ is a subobject of $A$, then $K(A)=\max\{K(B),K(A/B)\}$.
\item If there is an epimorphism $G^{(n)}\to A$ for some finite $n$, then $K(A)\leq K(G)$.
\item If $A\neq 0$ and $K(A)\neq\infty$, then $A$ has a critical subobject.
\item If $A$ is $\alpha$-critical, then any non-zero submodule of $A$ is also $\alpha$-critical.
\item If $A$ is critical, then $A$ is uniform.
\item If $R$ is a commutative noetherian ring, then $K(R)$ is equal to the classical Krull dimension (the maximal length of a chain of prime ideals).
\end{enumerate}
\end{prop}

Note the significance of part (3) in $\Mod\mhyphen R$, for $R$ a ring: it says that any finitely generated module has Krull dimension at most $K(R)$.

In this section, we investigate the relationship between Krull dimension and the geometric structure of the injective spectrum. The aim is to develop a picture where the Krull dimension of the ring is equal to the dimension of the spectrum (in the sense of the length of a maximal chain of irreducible closed subsets) and the dimension of the closure of a point $E$ in the spectrum is governed by the Krull dimension of $E$.

Unfortunately, this nice picture does not always work. We shall see it failing in Section \ref{quantum plane}. An open question is to establish conditions in which it does work.

Even before coming to such counterexamples to the nice picture we might hope for, a problem with the project proposed above is that injective objects, being `big', do not generally have a well-defined Krull dimension. However, we can circumvent this issue, at least in a locally noetherian Grothendieck category, by defining a closely related notion of dimension as follows.

For an arbitrary non-zero object $A$ in a locally noetherian Grothendieck category $\mathcal{A}$, let the {\bf critical dimension} $\cd(A)$ denote the minimum Krull dimension of non-zero subobjects of $A$. Note that, since every non-zero object has a finitely generated subobject and every finitely generated object is noetherian and so has Krull dimension (by part (1) of Proposition \ref{krull dim}), $A$ must contain a non-zero subobject whose Krull dimension is defined. Then, since Krull dimension is ordinal-valued, and so obeys the well-order property, $\cd(A)$ is well-defined for any non-zero $A$. The name will be explained shortly.

For the remainder of this section, $\mathcal{A}$ will always denote a locally noetherian Grothendieck category. We begin by noting two trivial but useful properties of critical dimension.

\begin{lemma}\label{crit and hom}
Let $A$ be any non-zero object of $\mathcal{A}$, and $B$ any object whose Krull dimension is defined. Then
\begin{enumerate}
\item $A$ has a $\cd(A)$-critical submodule;
\item If $K(B)<\cd(A)$, then $(B,A)=0$.
\end{enumerate}\end{lemma}

\proof
\begin{enumerate}
\item By definition of $\cd(A)$ and since $\mathcal{A}$ is locally noetherian, $A$ has a non-zero submodule $C$ with $K(C)=\cd(A)$. Any object which has Krull dimension contains a critical subobject, by part (4) of Proposition \ref{krull dim}, so there is $C'\leq C$ critical. Moreover, since $C'\leq A$, $K(C')\geq \cd(A)$, but since $C'\leq C$, $K(C')\leq K(C)=\cd(A)$; so $C'$ is $\cd(A)$-critical.
\item Suppose $K(B)<\cd(A)$ and let $\phi:B\to A$ be a morphism. Then $\phi(B)$ has Krull dimension at most $K(B)$, but $A$ has no non-zero subobject of Krull dimension less than $\cd(A)$, so $\phi(B)$ must be zero.\qed\end{enumerate}

\begin{lemma}\label{uniform critical sub}
Let $U$ be a uniform object of $\mathcal{A}$. Then any two critical subobjects of $U$ have the same Krull dimension.
\end{lemma}

\proof
Suppose $A,B\leq U$ are critical. Then, by Proposition \ref{krull dim}, any subobject of $A$ is $K(A)$-critical and any subobject of $B$ is $K(B)$-critical. But $U$ is uniform, so $A$ and $B$ have a common subobject, and therefore $K(A)=K(B)$.\qed

This critical dimension now allows us to apply Krull dimension to injective objects. By the above Lemma \ref{uniform critical sub}, for an indecomposable injective $E$, all critical subobjects must have the same Krull dimension; in light of Lemma \ref{crit and hom}, this dimension will be precisely $\cd(E)$. This is the reason for terming it the critical dimension.

This chapter, and our subsequent uses of the results herein, should be compared with results of Pappacena in \cite{pappacena}. There, it is shown that for any dimension function on a Grothendieck category satisfying suitable axioms, the sum of all critical subobjects of a uniform object is itself critical. Pappacena then works with the resulting largest critical subobject and relates its dimension to topological notions of dimension on the injective spectrum in his ``weak Zariski'' topology (see Section \ref{weak zariski}). Of course, given a uniform object $U$, the dimension of its largest critical subobject is precisely the critical dimension as defined here; so Pappacena uses essentially the same tool as we do, but reached via a slightly different route.

\subsection{Specialisation and Dimension; Points of Maximal Dimension; Irreducibility}

\begin{lemma}\label{specialisation alpha}
Let $E,F\in\injspec(R)$ be such that $E\leadsto F$. Then $\cd(E)\geq \cd(F)$, with equality if and only if $E\cong F$.\end{lemma}

\proof
Take a $\cd(E)$-critical subobject $C$ of $E$, by Lemma \ref{crit and hom}; then $E=E(C)$, since $E$ is indecomposable. Now $E\leadsto F$, so $(C,E)\neq 0$, so $(C,F)\neq 0$, by Lemma \ref{specialisation tfae}. It follows by Lemma \ref{crit and hom} that $\cd(E)=K(C)\geq \cd(F)$.

Take $\phi:C\to F$ non-zero. If $\phi$ is an embedding, then $F\cong E(C)=E$ (and certainly $\cd(E)=\cd(F)$). Otherwise, $\phi(C)$ is a proper quotient of $C$, hence
\[K(\phi(C))<K(C)=\cd(E),\]
by criticality of $C$. But $(\phi(C),F)\neq 0$, so, by Lemma \ref{crit and hom},
\[\cd(F)\leq K(\phi(C))<\cd(E).\]\qed

\begin{cor}\label{t0}
For any locally noetherian Grothendieck category $\mathcal{A}$, $\injspec(\mathcal{A})$ is $T_0$, \textit{i.e.}, Kolmogorov.
\end{cor}

\proof
If two points, $E$ and $F$, are topologically indistinguishable, then $E\leadsto F$ and $F\leadsto E$, so $\cd(E)\geq\cd(F)\geq\cd(E)$; hence $\cd(E)=\cd(F)$, but this yields an isomorphism $E\cong F$ when combined with the fact that $E\leadsto F$.\qed

We can obtain stronger results by strengthening the local noetherianity condition. Say that $\mathcal{A}$ is {\bf strongly locally noetherian via} $G$ if $\mathcal{A}$ is a Grothendieck category and $G$ is a generator for $\mathcal{A}$ which is noetherian. Note that a locally noetherian category has a generating \textit{set} of noetherian objects, but might not have any \textit{single} noetherian object which generates the whole category by itself; so not all locally noetherian categories are strongly locally noetherian. However, the category of right modules over a right noetherian ring $R$ is strongly locally noetherian via $R$.

Observe that, by part (3) of Proposition \ref{krull dim}, if $\mathcal{A}$ is strongly locally noetherian via $G_1$ and also via $G_2$, we must have $K(G_1)=K(G_2)$. So from the point of view of dimension, it does not matter which noetherian generator we take.

\begin{cor}\label{specialisation length}
Let $A$ be strongly locally noetherian via $G$, and write $d=K(G)$. Let $\alpha$ be an ordinal and let $E$ be a specialisation chain in $\injspec(\mathcal{A})$ of order type $\alpha^\op$; \textit{i.e.}, $E:\alpha\to\injspec(\mathcal{A})$ is an injection with $E_\gamma\leadsto E_\delta$ for all $\delta<\gamma<\alpha$. Then $\alpha\leq d+1$.
\end{cor}

\proof
By Lemma \ref{specialisation alpha}, for all $\delta<\gamma<\alpha$, $\cd(E_\gamma)>\cd(E_\delta)$. It follows by transfinite induction that $\cd(E_\gamma)\geq \gamma$. Moreover, since every indecomposable injective $E_\gamma$ contains a finitely generated object, whose Krull dimension is therefore at most $d$ (by part (3) of Proposition \ref{krull dim}), we see that $\cd(E_\gamma)\leq d$. So for each $\gamma<\alpha$, $d\geq \gamma$. If $\alpha$ is a limit ordinal, then it follows that $d\geq \alpha$.

If, on the other hand, $\alpha$ is a successor ordinal, say $\alpha=\gamma+1$, then $d\geq \cd(E_\gamma)\geq \gamma$, so $d+1\geq \alpha$.\qed

If $\injspec(\mathcal{A})$ is sober, then specialisation chains such as this correspond to chains of irreducible closed subsets. So for strongly locally noetherian categories where the injective spectrum is sober, the above corollary bounds the dimension of the space (in the sense of the length of a maximal chain of irreducible closed subsets) to be at most the dimension of a noetherian generator. See Section \ref{first examples} to see that noetherian rings of Krull dimension 0 (respectively 1) have injective spectra of dimension at most 0 (respectively 1).

Note that, since $\alpha$ in the above Corollary counts the number of injectives in a specialisation chain, not the number of specialisations, the bound on topological dimension really is $d$, not $d+1$. For instance, if $d=1$, then take $\alpha=2$, the maximum value allowed by the Corollary; then a chain as in the Corollary consists of indecomposable injectives $E_1\leadsto E_0$, which would normally be called a chain of length 1. So although $\alpha=d+1$, this does correspond to the maximum length of a specialisation chain being $d$.

Given this upper bound, it is natural to try to bound the dimension of the spectrum from below by $d$ as well. We now turn to this.

\begin{lemma}\label{dim -1 subquotient}
Fix an ordinal $\alpha$. Then:
\begin{enumerate}
\item If $A$ is an object of $\mathcal{A}$ with $K(A)=\alpha+1$, there is a subquotient of $A$ with Krull dimension exactly $\alpha$;
\item If $A$ is a non-zero, noetherian object with $K(A)=\alpha$, there is an $\alpha$-critical subquotient of $A$.
\end{enumerate}
\end{lemma}

\proof
\begin{enumerate}
\item Since $K(A)\geq \alpha+1$, there is an $\omega^\op$ descending chain $A_0>A_1>\hdots$ of subobjects of $A$ with $K(A_i/A_{i+1})\geq \alpha$ for all sufficiently large $i$. If for all but finitely many $i$ we had $K(A_i/A_{i+1})\geq \alpha+1$, then we would have $K(A)\geq \alpha+2$, a contradiction, so there is certainly some $i$ such that $K(A_i/A_{i+1})=\alpha$. Then this $A_i/A_{i+1}$ is a subquotient of $A$ with Krull dimension $\alpha$.

\item Let $A$ be a non-zero noetherian object of Krull dimension $\alpha$. Since $A$ is noetherian, the set of subobjects of $A$ with Krull dimension strictly less than $\alpha$ (which contains 0, so is non-empty) has a maximal element, $A_0$, say. By Proposition \ref{krull dim}, $A/A_0$ contains a critical subobject; we show that any non-zero subobject of $A/A_0$ has Krull dimension $\alpha$, and therefore this critical subobject is $\alpha$-critical, as required.

So let $B$ be a subobject of $A$ strictly containing $A_0$, so that $B/A_0$ is a general non-zero subobject of $A/A_0$. Then, by maximality of $A_0$, $K(B)=\alpha$; but by part (2) of Proposition \ref{krull dim}, $K(B)=\max\{K(A_0),K(B/A_0)\}$. Since $K(A_0)<\alpha$, we must have $K(B/A_0)=\alpha$, completing the proof.\qed\end{enumerate}

\begin{lemma}\label{subquotient support}
Let $A\leq B\leq C$ be a chain of objects in $\mathcal{A}$, so $B/A$ is a generic subquotient of $C$. Then for any injective object $E$, if $(B/A,E)\neq 0$, then $(C,E)\neq 0$. In particular, if $C$ is finitely presented, there is an inclusion of basic closed sets of $\injspec(\mathcal{A})$: $(B/A)\subseteq (C)$.\end{lemma}

\proof
Let $\phi:B/A\to E$ be a non-zero morphism. Then $\phi$ lifts to a non-zero map $B\to E$ and extends (by injectivity of $E$) to a non-zero map $C\to E$.\qed

\begin{lemma}\label{hull of critical}
Let $A$ be a critical object of $\mathcal{A}$. Then $E(A)$ is indecomposable and $\cd(E(A))=K(A)$.\end{lemma}

\proof
Critical objects are uniform by part (6) of Proposition \ref{krull dim}, so certainly $E(A)$ is indecomposable. Since $A\leq E(A)$, $\cd(E(A))\leq K(A)$. Moreover, if $B\leq E(A)$ has Krull dimension strictly less than $K(A)$, then $0\neq B\cap A\leq B$ and $K(A\cap B)\leq K(B)<K(A)$. But every non-zero submodule of $A$ has Krull dimension $K(A)$, by criticality, giving a contradiction. So $E(A)$ has no non-zero submodule with Krull dimension strictly less than $K(A)$, so $\cd(E(A))\geq K(A)$. The result follows.\qed

\begin{lemma}\label{direct sum basic open sets}
Let $\mathcal{A}$ be any Grothendieck category and
\[0\to A\to C\to B\to 0\]
a short exact sequence of finitely presented objects in $\mathcal{A}$. Then $(C)=(A)\cup (B)$.
\end{lemma}

\proof
Let $E$ be any injective object of $\mathcal{A}$. Then the functor $(-,E):\mathcal{A}^\op\to\ab$ is exact, giving us a short exact sequence of abelian groups:
\[0\to (B,E)\to (C,E)\to (A,E)\to 0.\]

The middle term in a short exact sequence is zero if and only if both outer terms are zero, so we see that $(C,E)=0$ if and only if $(A,E)=0$ and $(B,E)=0$. Since this holds for all indecomposable injectives $E$, this means $[C]=[A]\cap [B]$; taking complements gives the desired result.\qed

\begin{lemma}\label{irreducible implies critical}
Let $A$ be a non-zero, finitely presented object of $\mathcal{A}$ such that $(A)$ is an irreducible closed set. Then there is some finitely presented, critical object $B$ such that $(A)=(B)$.
\end{lemma}

\proof
First we let $S$ denote the set of critical subquotients of $A$ and show that
\[ (A)=\bigcup_{C\in S} (C).\]

The right-to-left inclusion is an immediate consequence of Lemma \ref{subquotient support}. For the reverse inclusion, suppose $E\in (A)$; then there is some non-zero map $\phi:A\to E$, and $\phi(A)$ has a critical subobject $C$, by part (4) of Proposition \ref{krull dim}. So $C\in S$ and $E\in (C)$, proving the equation.

Now we use the fact that $(A)$ is compact in the Ziegler topology (Lemma \ref{ziegler compact opens}), and hence it suffices to take the union over a finite subset of $S$. Irreducibility of $(A)$ then lets us narrow down further to a single element of $S$.\qed

So if $\mathcal{A}$ is strongly locally noetherian via $G$, with $d=K(G)<\omega$, the idea is to take a critical object $C_d$ of Krull dimension $d$, then apply Lemma \ref{dim -1 subquotient} to obtain a $(d-1)$-critical subquotient $C_{d-1}$. Continuing in this way, and applying Lemma \ref{subquotient support}, we obtain a chain of basic closed sets $(C_d)\supseteq (C_{d-1}) \supseteq \hdots \supseteq (C_0)$ of length $d$.

By Lemma \ref{hull of critical}, $(C_i)$ contains $E(C_i)$ and $\cd(E(C_i))=K(C_i)=i$; so we have $E(C_i)\notin (C_j)$ for any $j<i$, so the chain is proper. If therefore $(C_i)$ were irreducible for all $i$, this would bound the dimension of $\injspec(\mathcal{A})$ from below by $d$. Since Lemma \ref{irreducible implies critical} says that every irreducible basic closed set can be given by a critical module, we might hope that the converse holds, in which case each $C_i$ would indeed be irreducible. Unfortunately, we shall see in Section \ref{quantum plane} that it is possible for the basic closed set given by a critical module not to be irreducible, so the converse of Lemma \ref{irreducible implies critical} fails.

We can, however, extract the following result, giving a further link between Krull dimension and the Zariski topology on $\injspec(\mathcal{A})$.

\begin{prop}
Let $A$ and $B$ be finitely presented objects of $\mathcal{A}$ such that $(B)\subseteq(A)$. Then $K(B)\leq K(A)$.
\end{prop}

\proof
Since $B$ is noetherian, it has a $K(B)$-critical subquotient $B_0$, by Lemma \ref{dim -1 subquotient}. Then $E(B_0)\in (B_0)\subseteq (B)\subseteq (A)$, so $(A,E(B_0))\neq 0$, and so by Lemma \ref{crit and hom}, $K(A)\geq \cd(E(B_0))=K(B_0)=K(B)$.\qed

A noetherian object is artinian - \textit{i.e.} has Krull dimension 0 - if and only if it has a composition series, the factors of which are simple - \textit{i.e.}, 0-critical - objects. This generalises to higher dimensions, though with a weakening of the uniqueness.

For an object $A$ of $\mathcal{A}$, a {\bf critical composition series} is a finite chain $A=A_n>A_{n-1}>\hdots>A_1>A_0=0$ of subobjects such that for each $i$, $A_{i+1}/A_i$ is a critical object, and $K(A_{i+2}/A_{i+1})\geq K(A_{i+1}/A_i)$.

\begin{prop}[\cite{mcr}, 6.2.19-6.2.22]
Let $A$ be a noetherian object in $\mathcal{A}$. Then $A$ has a critical composition series. Moreover, any two critical composition series for $A$ have the same length, and their composition factors can be paired so that corresponding factors have a non-zero isomorphic subobject.
\end{prop}

Note that, since the composition factors are critical, and hence uniform, the uniqueness condition of the Proposition guarantees that the injective hulls of the composition factors are uniquely determined, even though the composition factors themselves need not be.

In a commutative noetherian ring, there are finitely many minimal prime ideals - \textit{i.e.}, finitely many irreducible closed sets in the Zariski spectrum of maximal dimension. We now show an analagous result for the injective spectrum of a strongly locally noetherian Grothendieck category.

\begin{thm}\label{finitely many maximal points}
Suppose $\mathcal{A}$ is strongly locally noetherian via $G$, with $K(G)=d$. Then there is at least one, and only finitely many, indecomposable injective objects of critical dimension $d$.
\end{thm}

\proof
By Lemma \ref{dim -1 subquotient}, there is a $d$-critical subquotient of $G$. The injective hull of this therefore has critical dimension $d$, so we have at least one such.

Now take a critical composition series for $G$:
\[G=G_n>\hdots>G_0=0\]
and let $E$ be an indecomposable injective of critical dimension $d$. First we show that for some $i$, $(G_{i+1}/G_i,E)\neq 0$. Certainly $(G_n,E)\neq 0$ and $(G_0,E)=0$, so there is some index $i$ such that $(G_{i+1},E)\neq 0$, but $(G_i,E)=0$. Then all maps $G_{i+1}\to E$ must vanish on $G_i$, hence induce maps on the quotient, so $(G_{i+1}/G_i,E)\neq 0$, as claimed.

Now let $\phi:G_{i+1}/G_i\to E$ be a non-zero map. Since $\cd(E)=d$, which is the maximum possible value for the Krull dimension of a finitely presented object, by part (3) of Proposition \ref{krull dim}, we must have $K(\phi(G_{i+1}/G_i))=d$, so $K(G_{i+1}/G_i)=d$. But $G_{i+1}/G_i$ is critical, so $\phi$ must be an embedding. Therefore $E=E(G_{i+1}/G_i)$. Since the critical composition series has finite length and uniquely determines the injective hulls of its composition factors, there can only be finitely many such $E$. Indeed, the number of such $E$ will be precisely the number of inequivalent factors in any critical composition series for $G$.\qed

There is also of course the question of when the injective spectrum is irreducible and, in this case, if it has a generic point. We have the following result in the case of modules over a ring.

\begin{thm}\label{domain has irred spec}
If $R$ is a right noetherian domain, then $\injspec(R)$ is irreducible and has $E(R_R)$ as a generic point.
\end{thm}

\proof
By classical results of Goldie, any right noetherian domain is right Ore, and so has uniform dimension 1; \textit{i.e.}, $R_R$ is uniform. So certainly $E(R_R)$ is indecomposable. We show that $E(R_R)$ is contained in every non-empty open set, which proves that it is a generic point and therefore the space is irreducible. It suffices to show that $E(R_R)$ is contained in every non-empty basic open set $[M]$ for $M$ finitely presented.

First we show that we can reduce to the case where $M$ is cyclic. Let $M$ be a finitely presented module such that $[M]\neq \varnothing$ and suppose that there is some $\phi:M\to E(R_R)$ non-zero. Then there is some $m\in M$ such that $\phi(m)\neq 0$, so $(mR,E(R_R))\neq 0$. Moreover, $[mR]\neq\varnothing$, since if $(mR,E)\neq 0$ for all $E\in\injspec(R)$, then, by injectivity, $(M,E)\neq 0$ for all $E$, so $[M]=\varnothing$, a contradiction. So if there is $M$ finitely presented with $E(R_R)\notin [M]\neq\varnothing$, then there is $mR$ cyclic such that $E(R_R)\notin [mR]\neq \varnothing$.

Now let $I$ be a proper right ideal of $R$, so that $[R/I]$ is non-empty, and suppose for a contradiction that $(R/I,E(R_R))\neq 0$. By enlarging $I$, we may assume without loss of generality that $R/I$ embeds in $E(R_R)$; write $f$ for such an embedding. Then $f(R/I)$ has non-zero intersection with $R_R$ in $E(R_R)$; so there are some $r\in R\take I$ and $s\in R\take0$ such that $f(r+I)$ and $s$ coincide in $E(R_R)$.

Then $\ann_R(r+I)=\ann_R(f(r+I))=\ann_R(s)=0$, since $R$ is a domain. So $(r+I)R\cong R_R$ is a free module of rank 1 inside $R/I$. We show that this cannot occur; \textit{i.e.}, that no proper quotient of $R_R$ contains an isomorphic copy of $R_R$.

Let $\phi_1:R_R\to R/I$ be an embedding. Then $(R/I)/(\phi_1(I))$ contains $\phi_1(R)/\phi_1(I)\cong R/I$. So let $I_1$ be the lift of $\phi_1(I)$ along the quotient map $R\to R/I$, so that $I_1$ strictly contains $I$ and $R/I_1$ contains an isomorphic copy of $R/I$. But since $R/I$ contains a copy of $R_R$, we have an embedding $\phi_2:R_R\to R/I_1$ and we can repeat this construction to obtain $I_2$ strictly containing $I_1$ with $R/I_2$ containing a copy of $R_R$.

Proceeding in this way, we construct an infinite, strictly ascending chain of right ideals $I<I_1<I_2<\hdots$, a contradiction. So we conclude that for any proper right ideal $I$, $(R/I,E(R_R))=0$; \textit{i.e.}, $E(R_R)\in [R/I]$.\qed

\begin{cor}\label{domain implies critical}
Let $R$ be a right noetherian domain. Then $R_R$ is a critical module.
\end{cor}

\proof
Suppose for a contradiction that $0<I<R_R$ is such that $K(R/I)=K(R)$. Then, by Lemma \ref{dim -1 subquotient}, $R/I$ has a $K(R)$-critical subquotient, say $J/K$ for some right ideals $J,K$, with $J>K\geq I>0$. Then $E(J/K)$ is indecomposable, so $E(R_R)\leadsto E(J/K)$, but $\cd(E(R_R))=K(R)=\cd(E(J/K))$, so $E(R_R)=E(J/K)$, by Lemma \ref{specialisation alpha}.

So $(J/K,E(R_R))\neq 0$, and hence $(R/K,E(R_R))\neq 0$. Therefore $[R/K]$ does not contain $E(R_R)$; by the Theorem, this forces $K=R_R$, but this is a contradiction, since $J>K$.\qed

\section{Examples}\label{first examples}

In this section, we develop some examples of the theory at its nicest, where the topology is well-behaved, Krull dimension corresponds nicely to topological dimension, and a more-or-less complete picture of the injective spectrum can be obtained. We also give a less tame example, where, for instance, the relationship between Krull dimension and the topology breaks down somewhat.

\subsection{Right Artinian Rings}\label{artinian}

This section addresses the simplest possible case; a ring $R$ is right artinian if and only if $K(R)=0$; so we consider 0-dimensional rings, which should, of course, be expected to have 0-dimensional spectra. This is indeed the case.

\begin{prop}
Let $R$ be right artinian. Then $\injspec(R)$ is a finite discrete space.
\end{prop}

\proof
Let $E$ be an indecomposable injective module. Then there is a finitely presented module $M$ such that $E=E(M)$. Then $M$ is artinian, so has a simple submodule $S$, so $E=E(S)$. Any simple module is annihilated by the Jacobson radical $J(R)$, and $R/J(R)$ is a finite direct sum of matrix rings over division rings, by Artin-Wedderburn, hence has finitely many simple modules. So there are only finitely many points in $\injspec(R)$.

Moreover, the injective hull of a simple module is a closed point, by Lemma \ref{socle implies closed}, and a finite space where every point is closed must be discrete.\qed

\subsection{1-Critical Rings}\label{1-critical}

Having dealt with 0-dimensional rings, in this section we consider rings $R$ such that $R_R$ is a 1-critical module. Recall Corollary \ref{domain implies critical}, which says that if $R$ is a right noetherian domain, then $R_R$ is critical. So any 1-dimensional, right noetherian domain is covered by the results of this section. We also have the following result, giving further examples.

\begin{lemma}[\cite{mcr}, 6.2.8]
Let $R$ be a hereditary noetherian prime ring, then $R$ is either artinian or is 1-critical.
\end{lemma}

\begin{thm}\label{spectrum of a 1-critical ring}
Let $R$ be a right noetherian ring such that $R_R$ is 1-critical. Then the indecomposable injective $R$-modules are $E(R_R)$ and the injective hulls of the simple modules. The open sets in $\injspec(R)$ are precisely the cofinite sets including $E(R_R)$, and of course the empty set. Therefore, $E(R_R)$ is a generic point in $\injspec(R)$.
\end{thm}

\proof
Since $R_R$ is 1-critical, $E(R_R)$ is indecomposable. All the other indecomposable injectives are $E(R/I)$ for some non-zero right ideal $I$; but each such $R/I$ is artinian, since $R_R$ is 1-critical, so each indecomposable injective except $E(R_R)$ is the hull of a simple module.

For the topology, we first show that the basic open sets $[M]$ for $M$ finitely presented are all either cofinite and include $E(R_R)$, or are empty. We proceed by induction on the number of generators in a generating set for $M$. Recall Lemma \ref{direct sum basic open sets}, which says that if $B$ is an extension of $A$ by $C$, then $[B]=[A]\cup [C]$.

The base case is when $M\cong R/I$ is cyclic. If $I=0$, $R/I=R$, which maps to everything, so $[R/I]$ is the empty set. If $I\neq 0$, then $R/I$ is artinian, since $R_R$ is 1-critical, so $R/I$ has a composition series, and the indecomposable injectives receiving a map from $R/I$ are precisely the injective hulls of the composition factors of $R/I$, by repeated applications of Lemma \ref{direct sum basic open sets}. So $[R/I]$ is cofinite. Moreover, since $K(R/I)<1=K(R)=\cd(E(R_R))$, we must have $(R/I,E(R_R))=0$, so $E(R_R)\in[R/I]$. This completes the base case.

Now suppose that the result is proved for $n$-generated modules and let $M$ be generated by $m_1,\hdots,m_{n+1}$. Define
\[N:=\sum_{i=1}^n m_iR,\]
so that we have a short exact sequence
\[0\to N\to M\to M/N\to 0\]
and $M/N$ is cyclic, being generated by $m_{n+1}+N$. Then the inductive hypothesis covers both $N$ and $M/N$, and Lemma \ref{direct sum basic open sets} again gives us that $[M]=[N]\cup[M/N]$ is cofinite and includes $E(R_R)$ or is empty, completing the induction.

Each open set is a union of basic open sets, hence also is cofinite and includes $E(R_R)$ (or is empty). So all open sets fit the description given. Finally, let $U\subseteq\injspec(R)$ be a cofinite set including $E(R_R)$; we show that $U$ is open, by showing that $C:=\injspec(R)\take U$ is closed.

Since $C$ is a finite set excluding $E(R_R)$, it contains only the injective hulls of some finitely many simple modules. By Lemma \ref{socle implies closed}, the hull of a simple is a closed point; since a finite union of closed sets is closed, we see that $C$ is indeed a closed set.\qed

In particular, let $k$ be a field of characteristic 0, and let $A_1(k)$ be the first Weyl algebra over $k$; that is:
\[A_1(k):=k\langle x,\partial\mid \partial x-x\partial = 1\rangle.\]
Then $A_1(k)$ is a noetherian domain of Krull dimension 1 (see \textit{e.g.}, \cite[{}6.6.8]{mcr}), so its injective spectrum fits the above description.

Note also that this description of the topology, having a collection of closed points and one generic point, is that of the affine line over any field. So the injective spectrum of a 1-critical ring is homeomorphic to an affine line as topological spaces, though there may not be any canonical way to identify it, even when there is a ``natural'' field to consider, and they will not generally by homeomorphic as ringed spaces. For instance, $\injspec(A_1(k))\cong \spec(k[x])$, but simply for the reason that they have the same cardinality (by Block's classification of the simple $A_1(k)$-modules \cite[Theorem 1]{block}) and the same, very simple topology - we also have $\injspec(A_1(k))\cong \spec(k(t)[x])$, for the same reason. It will sometimes (see Section \ref{heisenberg}) be convenient to picture the injective spectrum of a 1-critical ring as a line for the purposes of visualising the injective spectra of more complicated rings, but it should always be borne in mind that this is a fiction for ease of visualisation.

Observe that in both the 0-dimensional (artinian) and 1-critical cases, the injective spectrum is a noetherian topological space. We will see in subsection \ref{quantum plane} that this can fail in dimension 2.

\subsection{The Heisenberg Algebra}\label{heisenberg}

Let $k$ be an algebraically closed field of characteristic 0 and let $\mathfrak{h}$ denote the (first) Heisenberg algebra; \textit{viz}. the 3-dimensional Lie algebra with basis $\{p,q,z\}$, where $[p,q]=z,$ and $[z,-]=0$. Let $H$ be the universal enveloping algebra of $\mathfrak{h}$; so $H=k[z]\langle p,q\mid [p,q]=z\rangle.$ We will obtain a complete description of the points of $\injspec(H)$, and an extensive description of the topology.

For convenience, we let $H_\alpha:=H/(z-\alpha)H$ denote the quotient ring, for $\alpha\in k$, and let $f_\alpha:H\to H_\alpha$ denote the quotient map.

First we describe the points of $\injspec(H)$. Let $E$ be an indecomposable, injective $H$-module. There are two cases to consider; either there is some non-zero element of $k[z]$ which acts non-invertibly on $E$, or there isn't.

In the first case, suppose that $f(z)$ is such an element of minimal degree, and monic, without loss of generality. Because $E$ is injective (hence divisible, by Theorem \ref{equations in injectives} and the following remarks) and $H$ is a domain, $f(z)$ must act surjectively on $E$, so to be non-invertible, its action must be non-injective. So there is some $e\in E\take 0$ such that $ef(z)=0$.

Since $k$ is algebraically closed, we can factor $f(z)$ as a product of linear factors $(z-\alpha_i)$. Hitting $e$ with each $(z-\alpha_i)$ in turn, we find that there is a non-zero element of $E$ which is annihilated by one of the $(z-\alpha_i)$'s. Since $f$ was assumed to have minimal degree, we see that $f$ is linear.

So if any non-zero element of $k[z]$ acts non-invertibly on $E$, then there exists some $\alpha\in k$ such that $(z-\alpha)$ acts non-invertibly. Since $(z-\alpha)$ is central, the annihilator in $E$ of $(z-\alpha)$ is a submodule; so $E$ has a submodule $S$ which is a $H_\alpha$-module. So $E$ has the form $E(M_H)$ for some $M\in \Mod\mhyphen H_\alpha$. Since any essential extension of $M$ in $\Mod\mhyphen H_\alpha$ remains essential over $H$, we can replace $M$ by its injective hull over $H_\alpha$, which must be indecomposable, since $M$ is uniform. So $E$ has the form $E(F_H)$ for some $F\in\injspec(H_\alpha)$.

Thus, in the notation of Theorem \ref{induced by surjection}, those modules on which a non-zero element of $k[z]$ acts non-invertibly are of the form $f_\alpha^\ast(F)$ for some $\alpha\in k$ and $F\in\injspec(H_\alpha)$. By Theorem \ref{induced by surjection}, $f_\alpha^\ast$ is an embedding of topological spaces, so the internal topology of the set of points of $\injspec(H)$ containing an element annihilated by $(z-\alpha)$ is precisely the topology on $\injspec(H_\alpha)$.

Moreover, if $E$ contains an element annihilated by $z-\alpha$ and also an element annihilated by $z-\beta$, for some $\alpha,\beta\in k$, then the submodules $\ann_E(z-\alpha)$ and $\ann_E(z-\beta)$ intersect, by uniformity. So $E$ contains a non-zero element annihilated by $z-\alpha$ and $z-\beta$, and hence also by $\alpha-\beta$. Therefore we must have $\alpha=\beta$. So the set of indecomposable injectives of the form $f^\ast(F)$ for $F\in\injspec(H_\alpha)$ splits as a disjoint union over the different values of $\alpha\in k$.\gap

In the second case, every non-zero element of $k[z]$ acts invertibly on $E$. Let $S:=k[z]\take 0$; this is a central multiplicative (hence Ore) set in $H$; let $H_S$ denote the corresponding localised ring. Then $E$ is an indecomposable, injective $H_S$-module, so lies in the image of the continuous injection $\injspec(H_S)\to\injspec(H)$ induced by the (flat, epimorphic) localisation map of Corollary \ref{induced map}.

Since every point in $\injspec(H)$ must fall into (exactly) one of the above two cases, we see that, as a set, $\injspec(H)$ is the disjoint union of $\injspec(H_\alpha)$ as $\alpha$ ranges over $k$, and also $\injspec(H_S)$.\gap

We now describe each of these parts of the spectrum. Note that $H_\alpha$ has presentation $k\langle p,q\mid [p,q]=\alpha\rangle$. In the case where $\alpha\neq 0$, replacing $p$ by $p/\alpha$, we obtain the first Weyl algebra over $k$, whose injective spectrum was described in Theorem \ref{spectrum of a 1-critical ring} and the remarks immediately thereafter. If instead $\alpha=0$, $p$ and $q$ commute in the quotient, so we get $H_\alpha=k[p,q]$, whose spectrum is, of course, the affine plane.

We will henceforth refer to $\injspec(H_\alpha)$ (with $\alpha\neq 0$) and the corresponding subset of $\injspec(H)$ as being a line, for convenience; do not take this too literally however, as it is simply a topological space of cardinality $|k|$ whose closed sets are the finite sets omitting the generic. There is no canonical homeomorphism with the affine line on $k$, nor any canonical bijection between the closed points and any 1-dimensional vector space. It is simply convenient to visualise as linear.

Each of these subsets of $\injspec(H)$ for different $\alpha\in k$ is closed; for the image of $\injspec(H_\alpha)$ is precisely the basic closed set $(H/(z-\alpha)H)$ (including for $\alpha=0$). By Lemma \ref{induced by surjection}, the subspace topology on this closed set is precisely the (known) topology on $\injspec(H_\alpha)$.

Now consider $H_S$. This has presentation $k(z)\langle p,q\mid [p,q]=z\rangle$. Replacing $p$ by $p/z$, we get the first Weyl algebra over $k(z)$. By Corollary \ref{induced map}, the map $\injspec(H_S)\to\injspec(H)$ given by restricting scalars along the localisation map is an embedding of topological spaces. The image of this embedding, however, is not a closed set, as we shall shortly see.\gap

We now consider specialisation within $\injspec(H)$. Our tool here is Corollary \ref{specialisation lemma}; to show that $E\leadsto F$, we show that $E$ is torsionfree for $\mathcal{F}(F)$; for this, it suffices to prove that $M$ embeds in a direct product of copies of $F$, where $M$ is any module such that $E=E(M)$.

First we observe by Theorem \ref{domain has irred spec} that, since $H$ is a noetherian domain, $\injspec(H)$ is irreducible and $E(H)$ is generic in $\injspec(H)$. Observe that, as $H$ embeds in $H_S$, $E(H)=E(H_S)$, which is the generic point in the line $\injspec(H_S)$ embedded in $\injspec(H)$. In fact we shall see that the points in $\injspec(H_S)$ are generics over certain irreducible sets in the rest of the spectrum, and then $E(H_S)$ sits above all as the ``generic over the generics''.\gap

Let $\alpha\in k$ be non-zero. We show that for each closed point $E$ in $\injspec(H_\alpha)$, there is a point in $\injspec(H_S)$ (other than the generic) specialising to $E$. Let $I_\alpha$ be a maximal right ideal in $H_\alpha\cong A_1(k)$, so that $E(H_\alpha/I_\alpha)$ is a closed point over $H_\alpha$. Since $H_S\cong A_1(k(z))\cong A_1(k)\otimes_k k(z)$, $I_\alpha\otimes_k k(z)$ is a maximal proper right ideal of $H_S$ (after applying an isomorphism), so $E(H_S/(I_\alpha\otimes_k k(z)))$ is a closed point over $H_S$.

We show that $E(H_S/(I_\alpha\otimes_k k(z)))\leadsto E(H_\alpha/I_\alpha)$. Since $H$ embeds in $H_S$, we can consider the submodule
\[\frac{H+I_\alpha\otimes_k k(z)}{I_\alpha\otimes_k k(z)}\cong H/I,\]
where $I:=H\cap (I_\alpha\otimes_k k(z))$. Then
\[E\left(\frac{H_S}{I_\alpha\otimes_k k(z)}\right)=E(H/I).\]
So it suffices to prove that $H/I$ embeds in a direct product of copies of $E(H_\alpha/I_\alpha)$.

Let $\phi_n:H\to H/(I+(z-\alpha)^nH)$ be the quotient map. The intersection of the kernels of these $\phi_n$ over all $n\in\mathbb{N}$ is $I$, so if it can be shown that each $H/(I+(z-\alpha)^nH)$ embeds in $E(H_\alpha/I_\alpha)$ under some map $\psi_n$, then the product of the $\psi_n\circ \phi_n$'s will embed $H/I$ in $E(H_\alpha/I_\alpha)^{\aleph_0}$.

So we show that $H/(I+(z-\alpha)^n H)$ embeds in $E(H_\alpha/I_\alpha)$. This amounts to finding an element of $E(H_\alpha/I_\alpha)$ whose annihilator is exactly $I+(z-\alpha)^nH$. By solubility of equations in injective modules (Theorem \ref{equations in injectives}), there exists some $e\in E(H_\alpha/I_\alpha)$ such that $e(I+(z-\alpha)^nH)=0$ and $e(z-\alpha)^{n-1}=1+I_\alpha$ (which forces $e\neq 0$). Write $\xi_n(x)$ for the formula
\[x(I+(z-\alpha)^nH)=0\wedge x(z-\alpha)^{n-1}=1+I_\alpha.\]
We prove by induction on $n$ that if $E(H_\alpha/I_\alpha)\models\xi_n(e)$, then $\ann_H(e)= I+(z-\alpha)^nH$.

When $n=1$, $\xi_1(e)$ says that $e(I+(z-\alpha)H)=0$ and $e=1+I_\alpha$. This has a unique solution, $1+I_\alpha$, and $\ann_{H_\alpha}(1+I_\alpha)=I_\alpha$, so $\ann_H(1+I_\alpha)$ is the preimage of $I_\alpha$ under $f_\alpha$, which is $I+(z-\alpha)H$ (since a generating set for $I_\alpha$ lifts to a generating set for $I$), and so the base case is done.

For general $n\geq 2$, suppose any element satisfying $\xi_{n-1}$ has annihilator exactly $I+(z-\alpha)^{n-1}H$, and that $e$ satisfies $\xi_n$. Certainly $I+(z-\alpha)^nH\subseteq \ann_H(e)$, so we need only show the reverse inclusion. Suppose that $J$ is a right ideal containing $I+(z-\alpha)^nH$ and that $eJ=0$. We show that $J=I+(z-\alpha)^nH$.

In $H_\alpha$, $f_\alpha(I+(z-\alpha)H)=I_\alpha$, which is maximal, so $I+(z-\alpha)H$ is maximal in $H$, and hence $J+(z-\alpha)H=I+(z-\alpha)H$ or $J+(z-\alpha)H=H$.

If $J+(z-\alpha)H=H$, then there is some $h\in H$ such that $1-(z-\alpha)h\in J$. Then $(z-\alpha)^{n-1}-(z-\alpha)^nh\in J$, and
\[e[(z-\alpha)^{n-1}-(z-\alpha)^nh]=e(z-\alpha)^{n-1}-0=1+I_\alpha\neq 0,\]
so $J$ contains an element not annihilating $e$, a contradiction.

So $J+(z-\alpha)H=I+(z-\alpha)H$. Now, $e':=e(z-\alpha)$ satisfies
\[e'(I+(z-\alpha)^{n-1}H)=e((z-\alpha)I+(z-\alpha)^nH)=0\]
and
\[e'(z-\alpha)^{n-2}=e(z-\alpha)^{n-1}=1+I_\alpha,\]
so $e'$ satisifes $\xi_{n-1}$ and hence by the inductive hypothesis, $\ann_H(e')=I+(z-\alpha)^{n-1}H$.

For any $h\in J$, $h+(z-\alpha)g\in I$ for some $g\in H$, since $J+(z-\alpha)H=I+(z-\alpha)H$. Now, $eh=0$, but $e(z-\alpha)^{n-1}=1+I_\alpha$, so $e(z-\alpha)\neq 0$, as $n\geq 2$. But $eI=0$, so we must have $e(z-\alpha)g=0$; \textit{i.e.}, $e'g=0$. Therefore $g\in \ann_H(e')= I+(z-\alpha)^{n-1}H$. So
\[(z-\alpha)g\in I(z-\alpha)+(z-\alpha)^nH\subseteq I+(z-\alpha)^nH,\]
so $h=(h+(z-\alpha)g)-(z-\alpha)g\in I+(z-\alpha)^nH$, so $J=I+(z-\alpha)^nH$, completing the proof.

Note that, since each $H_\alpha$ for $\alpha\neq 0$ is isomorphic to $A_1(k)$, and we started with a maximal right ideal $I_\alpha$ in $H_\alpha$, we can take the corresponding right ideal in each $H_\beta$ for $\beta\neq 0$. So not only does $E(H/I)$ specialise to $E(H_\alpha/I_\alpha)$, but also to each corresponding point in each other fibre over $\beta\neq 0$. So, viewing each fibre as a vertical line, we can think of the closure of $E(H/I)$ as a horizontal line cutting across the fibres.\gap

So the picture we have developed of the injective spectrum of the universal enveloping algebra of the Heisenberg algebra is that we have, for each $\alpha\in k^\times$, a copy of the line $\injspec(A_1(k))$, together with an affine plane $\injspec(k[p,q])$ at $\alpha=0$, and a line $\injspec(A_1(k(z)))$ of generic points specialising across these $A_1(k)$-lines.

We can visualise this by imagining a copy of $k$ as a $z$-axis, with a fibre over each $\alpha\in k$, which is a line for non-zero $\alpha$, and a plane for $\alpha=0$. Each fibre is an irreducible closed set with its own generic point, and there is a ``line of generics'' which specialise across the fibres. The line of generics has its own ``big'' generic, $E(H)$, which specialises to every point. See Fig. 1.\gap

There are further points in the line of generics, with different closures cutting across the fibres. For instance, given any rational function $f(z)/g(z)\in k(z)$ (in lowest terms), the right ideal $(p-f(z)/g(z))H_S=(g(z)p-f(z))H_S$ is maximal in $H_S$, hence there is a point $E_{f/g}:=E(H/(g(z)p-f(z))H)$ in the line of generics. For any $\alpha$ such that $g(\alpha)\neq0$, we can take the $\cap$-irreducible right ideal $(g(\alpha)p-f(\alpha))H_\alpha$ of $H_\alpha$, and obtain the point $E(H_\alpha/(g(\alpha)p-f(\alpha))H_\alpha)=E(H/[(z-\alpha)H+(g(z)p-f(z))H])$, which we write as $E_{f/g,\alpha}$.

Then $E_{f/g}\leadsto E_{f/g,\alpha}$ for each $\alpha\in k$ such that $g(\alpha)\neq 0$. To prove this, we show that $H/(g(z)p-f(z))H$ embeds in a direct product of copies of $E_{f/g,\alpha}$, so that $E_{f/g}=E(H/(g(z)p-f(z))H)$ is torsionfree for $\mathcal{F}(E_{f/g,\alpha})$. It suffices to show that for each $n\geq 1$, $H/((z-\alpha)^nH+(g(z)p-f(z))H)$ embeds in $E_{f/g,\alpha}$; \textit{i.e.}, that $E_{f/g,\alpha}$ contains an element $e_n$ whose annihilator is precisely $(z-\alpha)^nH+(g(z)p-f(z))H$.

Let $\xi_n(x)$ be the system of equations
\[x(z-\alpha)^{n-1}=1+(z-\alpha)H+(g(\alpha)p-f(\alpha))H\wedge x(g(z)p-f(z))=0.\]
By solubility of systems of equations in injective modules (Theorem \ref{equations in injectives}), $\xi_n$ has a solution in $E_{f/g,\alpha}$, which must be non-zero, since $g(\alpha)\neq 0$ implies that $1+(z-\alpha)H+(g(\alpha)p-f(\alpha))H\neq 0$. We take $e_n$ to be any such solution and prove by induction that all solutions of $\xi_n$ have annihilator exactly $(z-\alpha)^nH+(g(z)p-f(z))H$.

The base case $n=1$ is clear, as the only solution to $\xi_1$ is $1+(z-\alpha)H+(g(\alpha)p-f(\alpha))H$, whose annihilator is exactly
\[(z-\alpha)H+(g(\alpha)p-f(\alpha))H=(z-\alpha)H+(g(z)p-f(z))H.\]
For $n\geq 2$, note that certainly $(z-\alpha)^nH+(g(z)p-f(z))H$ is contained in the annihilator of any solution $e_n$ of $\xi_n$, so we need only prove the reverse inclusion.

Suppose $h\in H$ is such that $e_nh=0$; we show that $h\in (z-\alpha)^nH+(g(z)p-f(z))H$. Certainly $e_n(z-\alpha)h=0$, and $e_n(z-\alpha)$ satisfies $\xi_{n-1}$; hence, by the inductive hypothesis, $h\in (z-\alpha)^{n-1}H+(g(z)p-f(z))H$. So we can write
\[h=(z-\alpha)^{n-1}h_1+(g(z)p-f(z))h_2\]
for some $h_1,h_2\in H$. Then
\[0=e_nh=e_n(z-\alpha)^{n-1}h_1+e_n(g(z)p-f(z))h_2=e_n(z-\alpha)^{n-1}h_1,\]
by $\xi_n$. But
\[e_n(z-\alpha)^{n-1}=1+(z-\alpha)H+(g(\alpha)p-f(\alpha))H,\]
and
\[(z-\alpha)H+(g(\alpha)p-f(\alpha))H=(z-\alpha)H+(g(z)p-f(z))H,\]
so this implies that $h_1=(z-\alpha)h_3+(g(z)p-f(z))h_4$, for some $h_3,h_4\in H$. So
\[h=(z-\alpha)^nh_3+(g(z)p-f(z))(h_2+(z-\alpha)^{n-1}h_4)\in (z-\alpha)^nH+(g(z)p-f(z))H,\]
as required.

So the point $E_{f/g}$ specialises to $E_{f/g,\alpha}$ whenever $g(\alpha)\neq 0$, as claimed. Note that when $g(\alpha)=0$, $(g(\alpha)p-f(\alpha))H=f(\alpha)H=H$ (since $f/g$ is assumed to be in lowest terms, so $f(\alpha)\neq 0$). So when $g(\alpha)=0$, $E_{f/g,\alpha}=0$, so is not a point in $\injspec(H)$.\gap

If $\alpha\neq 0$, then $H_\alpha\cong A_1(k)$ and, if $g(\alpha)\neq 0$, $(g(\alpha)p-f(\alpha))H_\alpha$ is a maximal ideal, so $E_{f/g,\alpha}$ is a closed point in the ``line'' $\injspec(H_\alpha)$. At $\alpha=0$, $H_0\cong k[p,q]$ and, if $g(0)\neq 0$, $(g(0)p-f(0))H_0$ is the prime ideal corresponding to the line $p=f(0)/g(0)$ in the affine plane $\injspec(H_0)\cong \spec(k[p,q])$. So $E_{f/g}$ specialises to the single closed point $E_{f/g,\alpha}$ for all $\alpha\neq 0$ such that $g(\alpha)\neq 0$, and to the generic point $E_{f/g,0}$ of the line $p=f(0)/g(0)$ (and hence to all the closed points of this line too) if $g(0)\neq 0$. Write $E_{f/g,0,\beta}$ for the closed point $E(H_0/((q-\beta)H_0+(g(0)p-f(0))H_0))$ in the closure of $E_{f/g,0}$, in the event that $g(0)\neq 0$.

Moreover, we now prove that these are the only points to which $E_{f/g}$ specialises, and in fact comprise the basic closed set $(H/(g(z)p-f(z))H)$. That is
\[\left(\frac{H}{(g(z)p-f(z))H}\right)=\cl(E_{f/g})=\{E_{f/g}\}\cup\{E_{f/g,\alpha}\mid g(\alpha)\neq 0\} \cup\{E_{f/g,0,\beta}\mid \beta\in k\},\]
where the final disjunct only occurs if $g(0)\neq 0$.

To see this, note that if $E$ is an indecomposable injective receiving a map from $H/(g(z)p-f(z))H$, then $E$ contains a non-zero element $e$ annihilated by $g(z)p-f(z)$. If $E$ lies in $\injspec(H_S)$, then $(g(z)p-f(z))H_S$ is a maximal right ideal of $H_S$, so $\ann_{H_S}(e)=(g(z)p-f(z))H_S$, and so $H_S/(g(z)p-f(z))H_S$ embeds in $E$, and hence $E=E_{f/g}$.

If instead $E$ lies in $\injspec(H_\alpha)$, then by the proof of Lemma \ref{induced by surjection} $E=E(F_H)$ for some indecomposable injective $H_\alpha$-module $F$ which receives a non-zero map from $H_\alpha/(g(z)p-f(z))H_\alpha=H_\alpha/(g(\alpha)p-f(\alpha))H_\alpha$. If $g(\alpha)=0$ then $f(\alpha)\neq 0$, since we assume that $f/g$ is in lowest terms, and so $H_\alpha/(g(\alpha)p-f(\alpha))H_\alpha=0$; so for those $\alpha$ such that $g(\alpha)=0$, there are no points of $(H/(g(z)p-f(z))H)$ lying in the fibre over $\alpha$.

So suppose $g(\alpha)\neq 0$. If $\alpha\neq 0$, $H_\alpha/(g(\alpha)p-f(\alpha))H_\alpha$ is simple, so the only possibility for $E$ is $E_{f/g,\alpha}$. If $\alpha=0$, $H_0\cong k[p,q]$ and $H_0/(g(0)p-f(0))H_0\cong k[p,q]/(p-f(0)/g(0))$, so $(H_0/(g(0)p-f(0))H_0)$ consists of the injective hulls of $k[p,q]/(p-f(0)/g(0))$ and $k[p,q]/(p-f(0)/g(0),q-\beta)$ for $\beta\in k$; but these points are precisely $E_{f/g,0}$ and $E_{f/g,0,\beta}$.

So every point in $(H/(g(z)p-f(z))H)$ is among the points listed, hence this is the entire set. Moreover, since $E_{f/g}\in (H/(g(z)p-f(z))H)$ and all the listed points are in $\cl(E_{f/g})$, we must also have that this set is equal to $\cl(E_{f/g})$.\gap

Similarly, we can consider $E(H_S/(g(z)q-f(z))H_S)$ and show that it specialises to each $E(H_\alpha/(g(\alpha)q-f(\alpha))H_\alpha)$ in each fibre for $\alpha\neq 0$ and $g(\alpha)\neq 0$ (and hence also to each closed point below this at $\alpha=0$), and that this is precisely the basic closed set $(H/(g(z)q-f(z))H)$.

Of course, not every maximal right ideal of $H_S\cong A_1(k(z))$ can be written in the form $(p-f(z)/g(z))A_1(k(z))$ for some rational function $f/g$ (or similarly with $q$ in lieu of $p$); so there are additional points in the line of generics, with perhaps more exotic closures. However, since each fibre is a closed set, the closure of a point can only pick out either finitely many points from each fibre over $\alpha\neq 0$, or the whole fibre. Over $\alpha=0$, the closure of a point can pick out a subvariety of the plane.

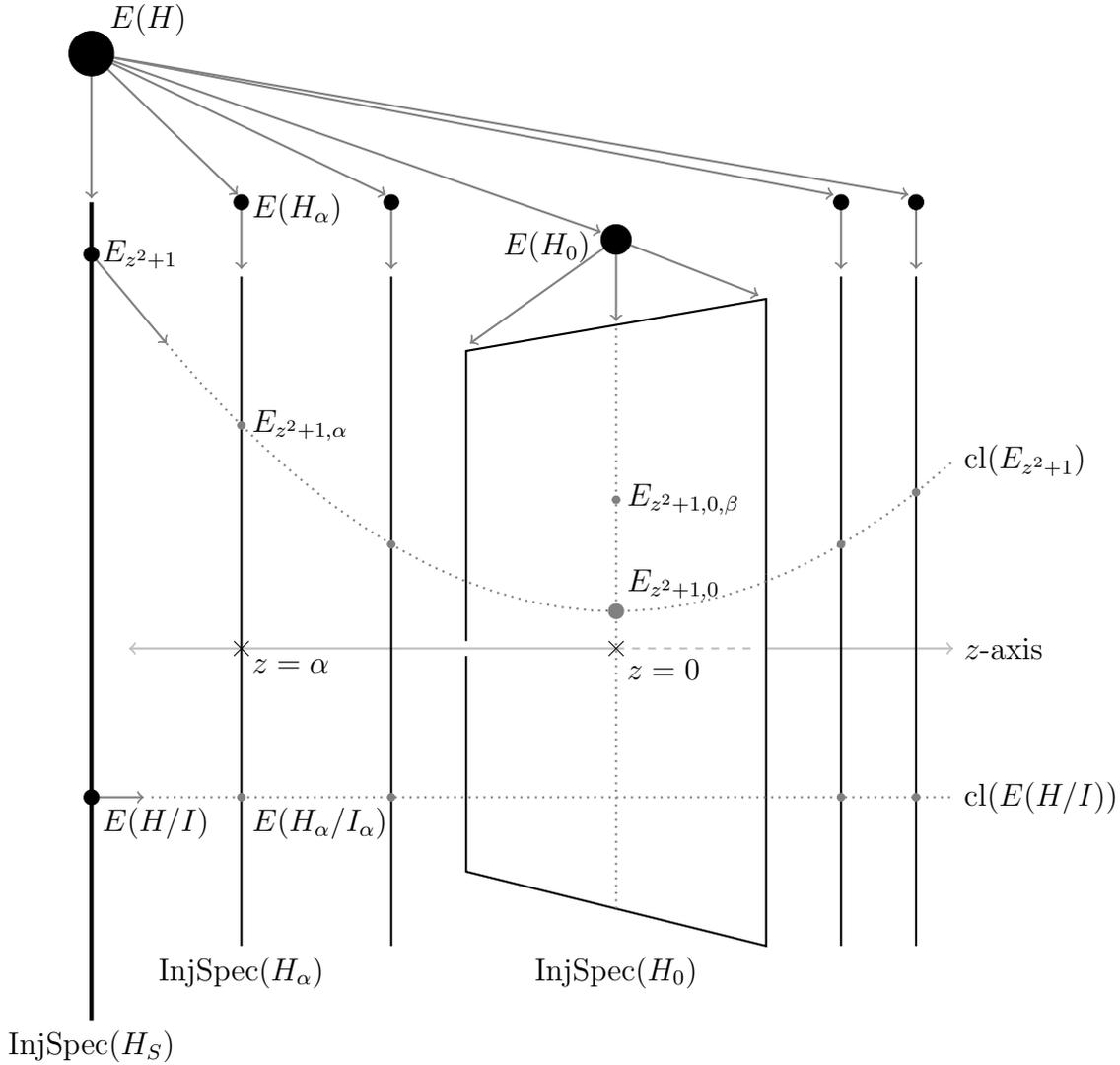
\begin{figure}
\begin{tikzpicture}

\draw[thick,lightgray,<-] (0.5,5) -- (7,5);
\draw[dashed, lightgray,thick] (7,5) -- (8.9,5);
\draw[thick,lightgray,->] (9,5) -- (11.5,5);
\node[right] at (11.5,5) {$z$-axis};

\draw[thick] (9,6) -- (9,1) -- (5,2) -- (5,4.9);
\draw[thick] (5,5.1) -- (5,9) -- (9,9.7) -- (9,6);
\node[below] at (7,1) {$\injspec(H_0)$};
\draw (6.9,5.1) -- (7.1,4.9); \draw (6.9,4.9) -- (7.1,5.1); 
\node[below right] at (7,5) {$z=0$};
\draw[thick,gray,->] (7,10.5) -- (5.05,9.1);
\draw[thick,gray,->] (7,10.5) -- (8.9,9.75);
\draw[thick,gray,->] (7,10.5) -- (7,9.4);
\draw[fill] (7,10.5) circle [radius=0.2];
\node[left] at (6.8,10.4) {$E(H_0)$};

\draw[thick] (2,1) -- (2,10);
\draw[thick] (4,1) -- (4,10);
\draw[thick] (10,1) -- (10,10);
\draw[thick] (11,1) -- (11,10);
\draw[thick,gray,->] (2,11) -- (2,10.1);
\draw[thick,gray,->] (4,11) -- (4,10.1);
\draw[thick,gray,->] (10,11) -- (10,10.1);
\draw[thick,gray,->] (11,11) -- (11,10.1);
\draw[fill] (2,11) circle [radius=0.1];
\draw[fill] (4,11) circle [radius=0.1];
\draw[fill] (10,11) circle [radius=0.1];
\draw[fill] (11,11) circle [radius=0.1];

\draw[ultra thick] (0,0) -- (0,11);
\node[below] at (0,0) {$\injspec(H_S)$};
\draw[->,gray,thick] (0,13) -- (0,11.05);
\draw[->,gray,thick] (0,13) -- (1.95,11.1);
\draw[->,gray,thick] (0,13) -- (3.93,11.1);
\draw[->,gray,thick] (0,13) -- (6.8,10.6);
\draw[->,gray,thick] (0,13) -- (9.9,11.1);
\draw[->,gray,thick] (0,13) -- (10.9,11.1);
\draw[fill] (0,13) circle [radius=0.3];
\node[above right] at (0.1,13.1) {$E(H)$};

\draw[thick,gray,dotted] (0.8,3) -- (11.5,3);
\node[right] at (11.5,3) {$\cl(E(H/I))$};
\draw[fill,gray] (2,3) circle [radius=0.05];
\draw[fill,gray] (4,3) circle [radius=0.05];
\draw[fill,gray] (10,3) circle [radius=0.05];
\draw[fill,gray] (11,3) circle [radius=0.05];
\draw[thick,gray,->] (0,3) -- (0.7,3);
\draw[fill] (0,3) circle [radius=0.1];
\node[below right] at (0,3) {$E(H/I)$};
\node[below right] at (2,3) {$E(H_\alpha/I_\alpha)$};

\draw[thick, gray, dotted, domain = 1:11.5] plot (\x, {5.5+0.1*(\x-7)*(\x-7)});
\draw[fill,gray] (2,8) circle [radius=0.05];
\node[right] at (2,8) {$E_{z^2+1,\alpha}$};
\draw[fill,gray] (4,6.4) circle [radius=0.05];
\draw[fill,gray] (10,6.4) circle [radius=0.05];
\draw[fill,gray] (11,7.1) circle [radius=0.05];
\draw[->,gray,thick] (0,10.3) -- (1,9.1);
\draw[fill] (0,10.3) circle [radius=0.1] node [right] {$E_{z^2+1}$};
\node[right] at (11.5,7.525) {$\cl(E_{z^2+1})$};
\draw[gray,dotted,thick] (7,1.5) -- (7,9.35);
\draw[fill,gray] (7,5.5) circle [radius=0.1];
\node[above right] at (7,5.5) {$E_{z^2+1,0}$};
\draw[fill,gray] (7,7) circle [radius=0.05];
\node[right] at (7,7) {$E_{z^2+1,0,\beta}$};

\node[below right] at (2,5) {$z=\alpha$};
\draw (1.9,5.1) -- (2.1,4.9); \draw (1.9,4.9) -- (2.1,5.1);
\node[below] at (2,1) {$\injspec(H_\alpha)$};
\node[right] at (2,10.9) {$E(H_\alpha)$};

\end{tikzpicture}
\caption{The injective spectrum of the Heisenberg algebra, shown as a collection of closed ``fibres'' at each value of $z$ in $k$, with a generic point over each fibre whose size indicates ``how generic it is'', and a ``line of generics'' to the left, with the ``biggest'' generic at the top. Specialisation is shown with grey arrows. The horizontal dotted grey line shows a point $E(H_\alpha/I_\alpha)$ at $z=\alpha\neq 0$ being lifted to a point $E(H/I)$ in the line of generics, which specialises to the original point and to the ``copies'' of it in the isomorphic fibres at other non-zero values of $z$. The dotted grey parabola and vertical line show the closure of a point $E_{z^2+1}$ in the line of generics to a single point from each fibre away from zero, and the line $p=0^2+1$ in $\injspec(H_0)\cong\spec(k[p,q])$; both the generic point of this line ($E_{z^2+1,0}$) and the closed point at $q=\beta$ ($E_{z^2+1,0,\beta}$) are shown.}
\end{figure}

\subsection{The Quantum Plane}\label{quantum plane}

Having seen examples of good behaviour from the injective spectrum, we now illustrate that the theory is not always so well-behaved. The following example is of a ring with Krull dimension 2 whose injective spectrum is not noetherian and where there is a closed point of critical dimension 1, and a critical module $M$ whose injective hull is not generic in $(M)$.

Let $k$ be an algebraically closed field and $q\in k$ a non-zero element that is not a root of unity. The quantum plane $A_q$ is the $k$-algebra
\[A_q:=k\langle x,y\mid xy=qyx\rangle,\]
which is the skew polynomial ring over $k[x]$ with new variable $y$, endomorphism defined by $x\mapsto q^{-1}x$, and derivation $0$. Note that $A_q$ is therefore noetherian, and has Krull dimension 2 (see \cite{mcr}, \S 6.9).

For any $\lambda\in k$, there is a simple $A_q$-module which we denote $k_{\lambda}$. This is a 1 dimensional $k$-vector space, where $x$ acts as scaling via $\lambda$ and $y$ via $0$. Each of these simple modules for different values of $\lambda$ is distinct; for two distinct eigenvalues of $x$ cannot occur in a one-dimensional space.

For any $\lambda\neq 0$, we will give an example of a finitely presented, 1-critical module $M_\lambda$ such that $E(k_{\lambda})\in (M_\lambda)$, but $E(M_\lambda)$ does not specialise to $E(k_{\lambda})$. This is a counterexample to the thought that for any finitely presented and critical module $M$, $E(M)$ is generic in $(M)$.

To construct $M_\lambda$, take a $k$-vector space with basis $\{v_i\mid i<\omega\}$ and define the action of $A_q$ by $v_ix=q^{-i}\lambda v_i$ and $v_iy=v_{i+1}$. Then $v_iyx=v_{i+1}x=q^{-i-1}\lambda v_{i+1}$, and $v_ixy=q^{-i}\lambda v_iy=q^{-i}\lambda v_{i+1}=q(v_iyx)$, so the commutation relation of $A_q$ is satisfied and this is indeed an $A_q$-module.

We see that $v_0$ is a generator for $M_\lambda$ as an $A_q$-module, so it is indeed finitely presented. Moreover, we will show that it is uniserial, with submodules $M_\lambda^{(n)}:=v_nA_q\cong M_{\lambda q^{-n}}$ for each $n<\omega$. To see this, suppose that $N$ is a submodule, containing some non-zero element
\[v=\sum_{i<\omega}\alpha_iv_i\]
for some $\alpha_i\in k$ almost all zero. For convenience, write $v=[\alpha_i]$.

Now $vx=[\alpha_i\lambda q^{-i}]$, so for any $n$ we have $v(x-q^{-n}\lambda)=[\alpha_i\lambda(q^{-i}-q^{-n})]$. But, since $\lambda\neq 0$ and $q$ is not a root of unity, $\alpha_i\lambda(q^{-i}-q^{-n})$ is zero if and only if $\alpha_i=0$ or $i=n$. So, taking $n$ such that $\alpha_n\neq 0$, $N$ contains an element $v(x-q^{-n}\lambda)$ whose representation in terms of the basis involves precisely one vector fewer than $v$ did. Repeating this process, we can eliminate all but one term from $v$ to obtain a non-zero element in $N$ which is a multiple of a single basis vector $v_n$ for some $n$. Moreover, we may do this for any $n$ such that $v_n$ is involved in a non-zero element of $N$; in particular, for the least such $n$.

So if $n$ is the least index such that $v_n$ is involved in an element of $N$, then $v_n\in N$. But then, by repeatedly applying $y$, we see that $v_i\in N$ for all $i\geq n$. So $N$ is precisely the $k$-linear span of $\{v_i\mid i\geq n\}$, or equivalently the cyclic submodule generated by $v_n$. Denote this cyclic submodule by $M_\lambda^{(n)}$. So the non-zero submodules of $M_\lambda$ are precisely the $M_\lambda^{(n)}$ for $n<\omega$.

Now, the $M_\lambda^{(n)}$ form an infinite descending chain of submodules of $M_\lambda$, so $K(M_\lambda)\geq 1$. Moreover, any quotient of $M_\lambda$ is $M_\lambda/M_\lambda^{(n)}$ for some $n$, which has $k$-dimension $n$, so is artinian. Therefore $M_\lambda$ is 1-critical.

The quotient $M_\lambda/M_\lambda^{(1)}$ is spanned by $\bar{v_0}=v_0+M_\lambda^{(1)}$, with $\bar{v_0}x=\lambda x$ and $\bar{v_0}y=0$, so this quotient is the simple module $k_{\lambda}$. So certainly $E(k_{\lambda})\in (M_\lambda)$.

We will now show that $E(M_\lambda)$ does not specialise to $E(k_{\lambda})$. Since $(M_\lambda^{(1)}, E(M_\lambda))\neq 0$, if $E(M_\lambda)\leadsto E(k_\lambda)$, then $(M_\lambda^{(1)},E(k_\lambda))\neq 0$, by Lemma \ref{specialisation tfae}. So we prove that $(M_\lambda^{(1)},E(k_\lambda))=0$.

Any proper, non-zero quotient of $M_\lambda^{(1)}$ has the form $M_\lambda^{(1)}/M_\lambda^{(n+1)}$ for some $n\geq 1$. But this contains as a submodule $M_\lambda^{(n)}/M_\lambda^{(n+1)}\cong k_{\lambda q^{-n}}$; so if $f:M_\lambda^{(1)}\to E(k_{\lambda})$ is non-zero, then either it is injective or its image contains the simple module $k_{\lambda q^{-n}}$ for some $n\geq 1$. But $E(k_\lambda)$ has $k_\lambda$ as an essential, simple submodule, so if $k_{\lambda q^{-n}}\in E(k_\lambda)$, then $k_{\lambda q^{-n}}\cong k_\lambda$, but this can only occur for $\lambda q^{-n}=\lambda$, contrary to the assumptions that $\lambda$ be non-zero and $q$ be not a root of unity.

So the only possibility for a non-zero map $M_\lambda^{(1)}\to E(k_\lambda)$ is an embedding. But then $M_\lambda^{(1)}$ must contain $k_\lambda$ as a submodule, whereas all submodules of $M_\lambda^{(1)}$ have the form $M_\lambda^{(n)}$ and so are infinite-dimensional.

Therefore we see that $(M_\lambda^{(1)},E(k_\lambda))=0$ and so $E(M_\lambda)$ does not specialise to $E(k_\lambda)$.

Furthermore, the point $E(M_\lambda)$, although of critical dimension 1, is a closed point, giving an example of a closed point which does not contain a simple submodule and thereby establishing that the converse to Lemma \ref{socle implies closed} is false, in general. For if $E(M_\lambda)\leadsto E\neq E(M_\lambda)$, then $\cd(E)<\cd(E(M_\lambda))$, so $\cd(E)=0$, and hence $E=E(S)$ for some simple module $S$. But then $(M_\lambda,E(S))\neq 0$ by Lemma \ref{specialisation tfae}, and so $S$ is a subquotient of $M_\lambda$. But the only simple subquotients of $M_\lambda$ are $k_{\lambda q^n}$ for $n<\omega$, and a straightforward adaptation of the above argument shows that $(M_\lambda^{(n+1)},E(k_{\lambda q^n}))=0$, so $E(M_\lambda)=E\left(M_\lambda^{(n+1)}\right)$ cannot specialise to $E(k_{\lambda q^n})$. So $E(M_\lambda)$ specialises to nothing (except of course itself).

Since the conjecture that any closed point is the hull of a simple module holds if $\injspec(R)$ is noetherian (Proposition \ref{basic closed implies socle}), this shows that $\injspec(A_q)$ cannot be noetherian. Indeed, the basic closed sets $\left(M_\lambda^{(n)}\right)$ form an infinite, strictly descending chain of closed sets, whose intersection is $\{E(M_\lambda)\}$.

\end{document}